\documentclass[a4paper,12pt]{article}
\usepackage{cmap}
\usepackage[T2A]{fontenc}
\usepackage[utf8]{inputenc} 
\usepackage[russian,english]{babel}
\usepackage{amsmath,amsthm,amssymb}
\usepackage{graphicx}
\usepackage{colortbl}
\usepackage{hyperref}
\usepackage{tikz}
\usepackage{currfile}
\usepackage{navigator}
\IfFileExists{\currfilename}{\embeddedfile{sourcefile}{\currfilename}}{}

\newtheorem{theorem}{Theorem}
\newtheorem{lemma}{Lemma}
\newtheorem{proposition}{Proposition}
\newtheorem{corollary}{Corollary}
\theoremstyle{remark}
\newtheorem{remark}{\bf Remark}

\makeatletter
\def\@seccntformat#1{\csname the#1\endcsname.\ } 
\def\@biblabel#1{#1.} 
\makeatother

\newcommand\vc[1]{\bar{#1}} 

\newcommand\quot[1]{\mbox{$\left(\begin{smallmatrix}#1\end{smallmatrix}\right)$}}
\newcommand\Quot[1]{\mbox{$\left(\begin{matrix}#1\end{matrix}\right)$}}
\newcommand\QQuot[1]{\mbox{$\Big(\begin{matrix}#1\end{matrix}\Big)$}}

\newcommand{\figureA}{%
\begin{figure}[ht]
\def\scaleX{2.0}
\def\scaleY{2.0}
\def\colora{green!20!blue!32}
\def\colorb{green!20!blue!32}
\def\1{$\boldsymbol 1$}
 \def\2{$\boldsymbol 1$}
 \def\aX{-0.4} \def\aY{+0.5}
 \def\bX{+0.0} \def\bY{+0.5}
 \def\cX{+0.4} \def\cY{+0.5}
 \def\dX{+1.2} \def\dY{+0.5}
 \def\eX{+2.2} \def\eY{+0.5}
 \def\XX{\aX\bX\cX\dX\eX-}
 \def\YY{\aY\bY\cY\dY\eY-}
 \def \bllsZ{15pt}
 \tikzstyle{vtX}=[circle,minimum size=12pt,inner sep=0pt]
 \tikzstyle{bll}=[circle,minimum size=\bllsZ,inner sep=0pt, fill=\colorb]
 \tikzstyle{bla}=[circle,minimum size=\bllsZ,inner sep=0pt, fill=\colora]
 \tikzstyle{blo}=[circle,minimum size=\bllsZ,inner sep=0pt, fill=white]
 \tikzstyle{vtx} = [vtX, thin, draw=black]
 \tikzstyle{vtC} = [vtX, fill=black]
 \tikzstyle{selected edge} = [draw,line width=5pt,-,red!50]
 \tikzstyle{edge} = [draw,-,black]
 \centering
 \scalebox{0.6}
{
 \begin{tikzpicture}[yscale=\scaleY,xscale=\scaleX]
\node[bla] at ( 0, 0 ) {};
\node[blo] at ( \aX, \aY) {};
\node[blo] at ( \aX\bX, \aY\bY) {};
\node[bla] (b11100) at ( \aX\bX\cX      , \aY\bY\cY      ) {};
\node[bla] (b11110) at ( \aX\bX\cX\dX   , \aY\bY\cY\dY   ) {};
\node[blo] (b01110) at (    \bX\cX\dX   ,    \bY\cY\dY   ) {};
\node[blo] at ( \cX\dX, \cY\dY) {};
\node[bla] at ( \dX, \dY) {};

\node[blo] at ( \bX, \bY) {};
\node[blo] at ( \cX, \cY) {};
\node[blo] at ( \aX\dX, \aY\dY) {};
\node[blo] at ( \bX\dX, \bY\dY) {};
\node[blo] at ( \aX\cX, \aY\cY) {};
\node[blo] at ( \bX\cX, \bY\cY) {};

 \node[blo] (b10110) at ( \aX   \cX\dX   , \aY   \cY\dY   ) {};
 \node[blo] (b11010) at ( \aX\bX   \dX   , \aY\bY   \dY   ) {};

 \node[vtx] (v00000) at ( 0              , 0               ) {\1}; 
 \node[vtx] (v00001) at ( \aX            , \aY             ) {}; 
 \node[vtx] (v00011) at ( \bX\aX         , \bY\aY          ) {}; 
 \node[vtx] (v00111) at ( \cX\bX\aX      , \cY\bY\aY       ) {\1}; 
 \node[vtx] (v01111) at ( \dX\cX\bX\aX   , \dY\cY\bY\aY    ) {\1}; 
 \node[vtx] (v01110) at ( \dX\cX\bX      , \dY\cY\bY       ) {}; 
 \node[vtx] (v01100) at ( \dX\cX         , \dY\cY          ) {}; 
 \node[vtx] (v01000) at ( \dX            , \dY             ) {\1}; 

 \node[vtx] (v00010) at ( \bX            , \bY             ) {}; 
 \node[vtx] (v00100) at ( \cX            , \cY             ) {}; 
 \node[vtx] (v00101) at ( \cX\aX         , \cY\aY          ) {}; 
 \node[vtx] (v00110) at ( \cX\bX         , \cY\bY          ) {}; 
 \node[vtx] (v01001) at ( \dX\aX         , \dY\aY          ) {}; 
 \node[vtx] (v01010) at ( \dX\bX         , \dY\bY          ) {}; 
 \node[vtx] (v01011) at ( \dX\bX\aX      , \dY\bY\aY       ) {}; 
 \node[vtx] (v01101) at ( \dX\cX\aX      , \dY\cY\aY       ) {}; 
 \draw[edge] (v00000) -- (v00010) -- (v00011) -- (v00001) -- (v00000);
 \draw[edge] (v00100) -- (v00110) -- (v00111) -- (v00101) -- (v00100);
 \draw[edge] (v01000) -- (v01010) -- (v01011) -- (v01001) -- (v01000);
 \draw[edge] (v01100) -- (v01110) -- (v01111) -- (v01101) -- (v01100);
 \draw[edge] (v00000) -- (v01000) -- (v01100) -- (v00100) -- (v00000);
 \draw[edge] (v00001) -- (v01001) -- (v01101) -- (v00101) -- (v00001);
 \draw[edge] (v00010) -- (v01010) -- (v01110) -- (v00110) -- (v00010);
 \draw[edge] (v00011) -- (v01011) -- (v01111) -- (v00111) -- (v00011);
 \end{tikzpicture}
 \qquad \begin{tikzpicture}[yscale=\scaleY,xscale=\scaleX]
\node[blo] at ( 0, 0 ) {};
\node[bla] at ( \aX, \aY) {};
\node[blo] at ( \aX\bX, \aY\bY) {};
\node[blo] (b11100) at ( \aX\bX\cX      , \aY\bY\cY      ) {};
\node[blo] (b11110) at ( \aX\bX\cX\dX   , \aY\bY\cY\dY   ) {};
\node[bla] (b01110) at (    \bX\cX\dX   ,    \bY\cY\dY   ) {};
\node[blo] at ( \cX\dX, \cY\dY) {};
\node[blo] at ( \dX, \dY) {};

\node[bla] at ( \bX, \bY) {};
\node[blo] at ( \cX, \cY) {};
\node[bla] at ( \aX\dX, \aY\dY) {};
\node[bla] at ( \bX\dX, \bY\dY) {};
\node[bla] at ( \aX\cX, \aY\cY) {};
\node[bla] at ( \bX\cX, \bY\cY) {};

 \node[bla] (b10110) at ( \aX   \cX\dX   , \aY   \cY\dY   ) {};
 \node[blo] (b11010) at ( \aX\bX   \dX   , \aY\bY   \dY   ) {};

 \node[vtx] (v00000) at ( 0              , 0               ) {}; 
 \node[vtx] (v00001) at ( \aX            , \aY             ) {\1}; 
 \node[vtx] (v00011) at ( \bX\aX         , \bY\aY          ) {}; 
 \node[vtx] (v00111) at ( \cX\bX\aX      , \cY\bY\aY       ) {}; 
 \node[vtx] (v01111) at ( \dX\cX\bX\aX   , \dY\cY\bY\aY    ) {}; 
 \node[vtx] (v01110) at ( \dX\cX\bX      , \dY\cY\bY       ) {\1}; 
 \node[vtx] (v01100) at ( \dX\cX         , \dY\cY          ) {}; 
 \node[vtx] (v01000) at ( \dX            , \dY             ) {}; 

 \node[vtx] (v00010) at ( \bX            , \bY             ) {\1}; 
 \node[vtx] (v00100) at ( \cX            , \cY             ) {}; 
 \node[vtx] (v00101) at ( \cX\aX         , \cY\aY          ) {\1}; 
 \node[vtx] (v00110) at ( \cX\bX         , \cY\bY          ) {\1}; 
 \node[vtx] (v01001) at ( \dX\aX         , \dY\aY          ) {\1}; 
 \node[vtx] (v01010) at ( \dX\bX         , \dY\bY          ) {\1}; 
 \node[vtx] (v01011) at ( \dX\bX\aX      , \dY\bY\aY       ) {}; 
 \node[vtx] (v01101) at ( \dX\cX\aX      , \dY\cY\aY       ) {\1}; 
 \draw[edge] (v00000) -- (v00010) -- (v00011) -- (v00001) -- (v00000);
 \draw[edge] (v00100) -- (v00110) -- (v00111) -- (v00101) -- (v00100);
 \draw[edge] (v01000) -- (v01010) -- (v01011) -- (v01001) -- (v01000);
 \draw[edge] (v01100) -- (v01110) -- (v01111) -- (v01101) -- (v01100);
 \draw[edge] (v00000) -- (v01000) -- (v01100) -- (v00100) -- (v00000);
 \draw[edge] (v00001) -- (v01001) -- (v01101) -- (v00101) -- (v00001);
 \draw[edge] (v00010) -- (v01010) -- (v01110) -- (v00110) -- (v00010);
 \draw[edge] (v00011) -- (v01011) -- (v01111) -- (v00111) -- (v00011);
 \end{tikzpicture}
 \qquad
 \begin{tikzpicture}[yscale=\scaleY,xscale=\scaleX]
\node[bll] at ( 0, 0 ) {};
\node[bll] at ( \aX, \aY) {};
\node[bla] at ( \aX\bX, \aY\bY) {};
\node[bla] (b11100) at ( \aX\bX\cX      , \aY\bY\cY      ) {};
\node[bll] (b11110) at ( \aX\bX\cX\dX   , \aY\bY\cY\dY   ) {};
\node[bll] (b01110) at (    \bX\cX\dX   ,    \bY\cY\dY   ) {};
\node[bla] at ( \cX\dX, \cY\dY) {};
\node[bla] at ( \dX, \dY) {};

\node[blo] at ( \bX, \bY) {};
\node[blo] at ( \cX, \cY) {};
\node[blo] at ( \aX\dX, \aY\dY) {};
\node[blo] at ( \bX\dX, \bY\dY) {};
\node[blo] at ( \aX\cX, \aY\cY) {};
\node[blo] at ( \bX\cX, \bY\cY) {};
 
 \node[blo] (b10110) at ( \aX   \cX\dX   , \aY   \cY\dY   ) {};
 \node[blo] (b11010) at ( \aX\bX   \dX   , \aY\bY   \dY   ) {};
 
 \node[vtx] (v00000) at ( 0              , 0               ) {\2}; 
 \node[vtx] (v00001) at ( \aX            , \aY             ) {\2}; 
 \node[vtx] (v00011) at ( \bX\aX         , \bY\aY          ) {\1}; 
 \node[vtx] (v00111) at ( \cX\bX\aX      , \cY\bY\aY       ) {\1}; 
 \node[vtx] (v01111) at ( \dX\cX\bX\aX   , \dY\cY\bY\aY    ) {\2}; 
 \node[vtx] (v01110) at ( \dX\cX\bX      , \dY\cY\bY       ) {\2}; 
 \node[vtx] (v01100) at ( \dX\cX         , \dY\cY          ) {\1}; 
 \node[vtx] (v01000) at ( \dX            , \dY             ) {\1}; 
 
 \node[vtx] (v00010) at ( \bX            , \bY             ) {}; 
 \node[vtx] (v00100) at ( \cX            , \cY             ) {}; 
 \node[vtx] (v00101) at ( \cX\aX         , \cY\aY          ) {}; 
 \node[vtx] (v00110) at ( \cX\bX         , \cY\bY          ) {}; 
 \node[vtx] (v01001) at ( \dX\aX         , \dY\aY          ) {}; 
 \node[vtx] (v01010) at ( \dX\bX         , \dY\bY          ) {}; 
 \node[vtx] (v01011) at ( \dX\bX\aX      , \dY\bY\aY       ) {}; 
 \node[vtx] (v01101) at ( \dX\cX\aX      , \dY\cY\aY       ) {}; 
 \draw[edge] (v00000) -- (v00010) -- (v00011) -- (v00001) -- (v00000);
 \draw[edge] (v00100) -- (v00110) -- (v00111) -- (v00101) -- (v00100);
 \draw[edge] (v01000) -- (v01010) -- (v01011) -- (v01001) -- (v01000);
 \draw[edge] (v01100) -- (v01110) -- (v01111) -- (v01101) -- (v01100);
 \draw[edge] (v00000) -- (v01000) -- (v01100) -- (v00100) -- (v00000);
 \draw[edge] (v00001) -- (v01001) -- (v01101) -- (v00101) -- (v00001);
 \draw[edge] (v00010) -- (v01010) -- (v01110) -- (v00110) -- (v00010);
 \draw[edge] (v00011) -- (v01011) -- (v01111) -- (v00111) -- (v00011);
 \end{tikzpicture}
 }
 \caption{The three nonequivalent degree-$2$  perfect $2$-colorings of $Q_4$}
 \label{f:F0} 
 \end{figure}
} 

\newcommand{\figureB}{%
\begin{figure}[ht]
 \def\1{$\boldsymbol 1$}
 \def\2{$\boldsymbol 2$}
 \def\3{$\boldsymbol 3$}
 \def\4{$\boldsymbol 4$}
\def\colora{green!70!black!40}
\def\colorb{red!32}
\def\colorc{red!40!yellow!42}
\def\colord{yellow!52}
\def\colore{cyan!70!blue!37}
\def\colorf{violet!70!blue!27}
\def\styleDefs{%
 \tikzstyle{vtX}=[circle,minimum size=12pt,inner sep=0pt]
 \tikzstyle{bll}=[circle,minimum size=\bllsZ,inner sep=0pt, fill=\colorb]
 \tikzstyle{bla}=[circle,minimum size=\bllsZ,inner sep=0pt, fill=\colora]
 \tikzstyle{blb}=[circle,minimum size=\bllsZ,inner sep=0pt, fill=\colorc]
 \tikzstyle{blo}=[circle,minimum size=\bllsZ,inner sep=0pt, fill=\colord]
 \tikzstyle{bls}=[circle,minimum size=\bllsZ,inner sep=0pt, fill=\colore]
 \tikzstyle{blt}=[circle,minimum size=\bllsZ,inner sep=0pt, fill=\colorf]
 \tikzstyle{vtx} = [vtX, thin, draw=black]
 \tikzstyle{vtC} = [vtX, fill=black]
 \tikzstyle{selected edge} = [draw,line width=5pt,-,red!50]
 \tikzstyle{edge} = [draw,-,black]%
}
 \def\aX{-0.4} \def\aY{+0.5}
 \def\bX{+0.0} \def\bY{+0.5}
 \def\cX{+0.4} \def\cY{+0.5}
 \def\dX{+1.1} \def\dY{+0.5}
 \def\eX{+2.2} \def\eY{+0.5}
 \def\XX{\aX\bX\cX\dX\eX-}
 \def\YY{\aY\bY\cY\dY\eY-}
 \def \bllsZ{15pt}
 \centering
 \scalebox{0.6}
 {
 \begin{tikzpicture}[yscale=2.5,xscale=2.0]
\styleDefs
\node[bla] at ( 0, 0 ) {};
\node[bla] at ( \aX, \aY) {};
\node[blo] at ( \aX\bX, \aY\bY) {};
\node[bll] (b11100) at ( \aX\bX\cX      , \aY\bY\cY      ) {};
\node[bla] (b11110) at ( \aX\bX\cX\dX   , \aY\bY\cY\dY   ) {};
\node[bla] (b01110) at (    \bX\cX\dX   ,    \bY\cY\dY   ) {};
\node[blo] at ( \cX\dX, \cY\dY) {};
\node[bll] at ( \dX, \dY) {};

\node[blo] at ( \bX, \bY) {};
\node[blb] at ( \cX, \cY) {};
\node[bll] at ( \aX\dX, \aY\dY) {};
\node[blb] at ( \bX\dX, \bY\dY) {};
\node[blb] at ( \aX\cX, \aY\cY) {};
\node[bll] at ( \bX\cX, \bY\cY) {};
 
 \node[blo] (b10110) at ( \aX   \cX\dX   , \aY   \cY\dY   ) {};
 \node[blb] (b11010) at ( \aX\bX   \dX   , \aY\bY   \dY   ) {};

 \node[vtx,label=left:$_{0000}$] (v00000) at ( 0              , 0               ) {\1}; 
 \node[vtx,label=below:$_{1000}$] (v00001) at ( \aX , \aY  ) {\1}; 
 \node[vtx] (v00011) at ( \bX\aX         , \bY\aY          ) {\2}; 
 \node[vtx,label=above:$_{1110}$] (v00111) at ( \cX\bX\aX, \cY\bY\aY) {\4}; 
 \node[vtx] (v01111) at ( \dX\cX\bX\aX   , \dY\cY\bY\aY    ) {\1}; 
 \node[vtx] (v01110) at ( \dX\cX\bX      , \dY\cY\bY       ) {\1}; 
 \node[vtx] (v01100) at ( \dX\cX         , \dY\cY          ) {\2}; 
 \node[vtx,label=below:$_{0001}$] (v01000) at ( \dX , \dY  ) {\4}; 

 \node[vtx] (v00010) at ( \bX            , \bY             ) {\2}; 
 \node[vtx,label=below:$_{0010}$] (v00100) at ( \cX , \cY  ) {\3}; 
 \node[vtx] (v00101) at ( \cX\aX         , \cY\aY          ) {\3}; 
 \node[vtx] (v00110) at ( \cX\bX         , \cY\bY          ) {\4}; 
 \node[vtx] (v01001) at ( \dX\aX         , \dY\aY          ) {\4}; 
 \node[vtx] (v01010) at ( \dX\bX         , \dY\bY          ) {\3}; 
 \node[vtx] (v01011) at ( \dX\bX\aX      , \dY\bY\aY       ) {\3}; 
 \node[vtx] (v01101) at ( \dX\cX\aX      , \dY\cY\aY       ) {\2}; 
 \draw[edge] (v00000) -- (v00010) -- (v00011) -- (v00001) -- (v00000);
 \draw[edge] (v00100) -- (v00110) -- (v00111) -- (v00101) -- (v00100);
 \draw[edge] (v01000) -- (v01010) -- (v01011) -- (v01001) -- (v01000);
 \draw[edge] (v01100) -- (v01110) -- (v01111) -- (v01101) -- (v01100);
 \draw[edge] (v00000) -- (v01000) -- (v01100) -- (v00100) -- (v00000);
 \draw[edge] (v00001) -- (v01001) -- (v01101) -- (v00101) -- (v00001);
 \draw[edge] (v00010) -- (v01010) -- (v01110) -- (v00110) -- (v00010);
 \draw[edge] (v00011) -- (v01011) -- (v01111) -- (v00111) -- (v00011);
 \end{tikzpicture}
\qquad
 \begin{tikzpicture}[yscale=2.5,xscale=2.0]
\styleDefs
\node[bls] at ( 0, 0 ) {};
\node[bls] at ( \aX, \aY) {};
\node[bla] at ( \aX\bX, \aY\bY) {};
\node[bla] (b11100) at ( \aX\bX\cX      , \aY\bY\cY      ) {};
\node[bls] (b11110) at ( \aX\bX\cX\dX   , \aY\bY\cY\dY   ) {};
\node[bls] (b01110) at (    \bX\cX\dX   ,    \bY\cY\dY   ) {};
\node[bla] at ( \cX\dX, \cY\dY) {};
\node[bla] at ( \dX, \dY) {};
\node[blo] at ( \bX, \bY) {};
\node[blb] at ( \cX, \cY) {};
\node[blo] at ( \aX\dX, \aY\dY) {};
\node[blb] at ( \bX\dX, \bY\dY) {};
\node[blb] at ( \aX\cX, \aY\cY) {};
\node[blo] at ( \bX\cX, \bY\cY) {};
 \node[blo] (b10110) at ( \aX   \cX\dX   , \aY   \cY\dY   ) {};
 \node[blb] (b11010) at ( \aX\bX   \dX   , \aY\bY   \dY   ) {};

 \node[vtx] (v00000) at ( 0              , 0               ) {\1}; 
 \node[vtx] (v00001) at ( \aX            , \aY             ) {\1}; 
 \node[vtx] (v00011) at ( \bX\aX         , \bY\aY          ) {\2}; 
 \node[vtx] (v00111) at ( \cX\bX\aX      , \cY\bY\aY       ) {\2}; 
 \node[vtx] (v01111) at ( \dX\cX\bX\aX   , \dY\cY\bY\aY    ) {\1}; 
 \node[vtx] (v01110) at ( \dX\cX\bX      , \dY\cY\bY       ) {\1}; 
 \node[vtx] (v01100) at ( \dX\cX         , \dY\cY          ) {\2}; 
 \node[vtx] (v01000) at ( \dX            , \dY             ) {\2}; 

 \node[vtx] (v00010) at ( \bX            , \bY             ) {\3}; 
 \node[vtx] (v00100) at ( \cX            , \cY             ) {\4}; 
 \node[vtx] (v00101) at ( \cX\aX         , \cY\aY          ) {\4}; 
 \node[vtx] (v00110) at ( \cX\bX         , \cY\bY          ) {\3}; 
 \node[vtx] (v01001) at ( \dX\aX         , \dY\aY          ) {\3}; 
 \node[vtx] (v01010) at ( \dX\bX         , \dY\bY          ) {\4}; 
 \node[vtx] (v01011) at ( \dX\bX\aX      , \dY\bY\aY       ) {\4}; 
 \node[vtx] (v01101) at ( \dX\cX\aX      , \dY\cY\aY       ) {\3}; 
 \draw[edge] (v00000) -- (v00010) -- (v00011) -- (v00001) -- (v00000);
 \draw[edge] (v00100) -- (v00110) -- (v00111) -- (v00101) -- (v00100);
 \draw[edge] (v01000) -- (v01010) -- (v01011) -- (v01001) -- (v01000);
 \draw[edge] (v01100) -- (v01110) -- (v01111) -- (v01101) -- (v01100);
 \draw[edge] (v00000) -- (v01000) -- (v01100) -- (v00100) -- (v00000);
 \draw[edge] (v00001) -- (v01001) -- (v01101) -- (v00101) -- (v00001);
 \draw[edge] (v00010) -- (v01010) -- (v01110) -- (v00110) -- (v00010);
 \draw[edge] (v00011) -- (v01011) -- (v01111) -- (v00111) -- (v00011);
 \end{tikzpicture}
 \qquad
 \begin{tikzpicture}[yscale=2.5,xscale=2.0]
\styleDefs
\node[blt] at ( 0, 0 ) {};
\node[blt] at ( \aX, \aY) {};
\node[bls] at ( \aX\bX, \aY\bY) {};
\node[bls] (b11100) at ( \aX\bX\cX      , \aY\bY\cY      ) {};
\node[blt] (b11110) at ( \aX\bX\cX\dX   , \aY\bY\cY\dY   ) {};
\node[blt] (b01110) at (    \bX\cX\dX   ,    \bY\cY\dY   ) {};
\node[bls] at ( \cX\dX, \cY\dY) {};
\node[bls] at ( \dX, \dY) {};
\node[blo] at ( \bX, \bY) {};
\node[bla] at ( \cX, \cY) {};
\node[bla] at ( \aX\dX, \aY\dY) {};
\node[blo] at ( \bX\dX, \bY\dY) {};
\node[blo] at ( \aX\cX, \aY\cY) {};
\node[bla] at ( \bX\cX, \bY\cY) {};
 \node[blo] (b10110) at ( \aX   \cX\dX   , \aY   \cY\dY   ) {};
 \node[bla] (b11010) at ( \aX\bX   \dX   , \aY\bY   \dY   ) {};

 \node[vtx] (v00000) at ( 0              , 0               ) {\1}; 
 \node[vtx] (v00001) at ( \aX            , \aY             ) {\1}; 
 \node[vtx] (v00011) at ( \bX\aX         , \bY\aY          ) {\2}; 
 \node[vtx] (v00111) at ( \cX\bX\aX      , \cY\bY\aY       ) {\2}; 
 \node[vtx] (v01111) at ( \dX\cX\bX\aX   , \dY\cY\bY\aY    ) {\1}; 
 \node[vtx] (v01110) at ( \dX\cX\bX      , \dY\cY\bY       ) {\1}; 
 \node[vtx] (v01100) at ( \dX\cX         , \dY\cY          ) {\2}; 
 \node[vtx] (v01000) at ( \dX            , \dY             ) {\2}; 
 
 \node[vtx] (v00010) at ( \bX            , \bY             ) {\4}; 
 \node[vtx] (v00100) at ( \cX            , \cY             ) {\3}; 
 \node[vtx] (v00101) at ( \cX\aX         , \cY\aY          ) {\4}; 
 \node[vtx] (v00110) at ( \cX\bX         , \cY\bY          ) {\3}; 
 \node[vtx] (v01001) at ( \dX\aX         , \dY\aY          ) {\3}; 
 \node[vtx] (v01010) at ( \dX\bX         , \dY\bY          ) {\4}; 
 \node[vtx] (v01011) at ( \dX\bX\aX      , \dY\bY\aY       ) {\3}; 
 \node[vtx] (v01101) at ( \dX\cX\aX      , \dY\cY\aY       ) {\4}; 
 \draw[edge] (v00000) -- (v00010) -- (v00011) -- (v00001) -- (v00000);
 \draw[edge] (v00100) -- (v00110) -- (v00111) -- (v00101) -- (v00100);
 \draw[edge] (v01000) -- (v01010) -- (v01011) -- (v01001) -- (v01000);
 \draw[edge] (v01100) -- (v01110) -- (v01111) -- (v01101) -- (v01100);
 \draw[edge] (v00000) -- (v01000) -- (v01100) -- (v00100) -- (v00000);
 \draw[edge] (v00001) -- (v01001) -- (v01101) -- (v00101) -- (v00001);
 \draw[edge] (v00010) -- (v01010) -- (v01110) -- (v00110) -- (v00010);
 \draw[edge] (v00011) -- (v01011) -- (v01111) -- (v00111) -- (v00011);
 \end{tikzpicture}
 }
 \caption{The perfect $4$-colorings of $Q_4$ with quotient matrix $[[1,1,1,1],[1,1,1,1],[1,1,1,1],[1,1,1,1]]$}
 \label{f:H42} 
 \end{figure}%
} 

\newcommand{\figureC}{%
\begin{figure}[ht]
\def\colora{green!20!blue!32}
\def\colorb{yellow!64}
\def\1{$\boldsymbol 1$}
 \def\2{$\boldsymbol 2$}
 \def\aX{-0.4} \def\aY{+0.5}
 \def\bX{+0.0} \def\bY{+0.5}
 \def\cX{+0.4} \def\cY{+0.5}
 \def\dX{+1.1} \def\dY{+0.5}
 \def\eX{+2.2} \def\eY{+0.5}
 \def\XX{\aX\bX\cX\dX\eX-}
 \def\YY{\aY\bY\cY\dY\eY-}
 \def \bllsZ{15pt}
 \tikzstyle{vtX}=[circle,minimum size=12pt,inner sep=0pt]
 \tikzstyle{bll}=[circle,minimum size=\bllsZ,inner sep=0pt, fill=\colorb]
 \tikzstyle{bla}=[circle,minimum size=\bllsZ,inner sep=0pt, fill=\colora]
 \tikzstyle{blo}=[circle,minimum size=\bllsZ,inner sep=0pt, fill=white]
 \tikzstyle{vtx} = [vtX, thin, draw=black]
 \tikzstyle{vtC} = [vtX, fill=black]
 \tikzstyle{selected edge} = [draw,line width=5pt,-,red!50]
 \tikzstyle{edge} = [draw,-,black]
 \centering
 \scalebox{0.6}
{
 \begin{tikzpicture}[yscale=2.5,xscale=2.0]
\node[bla] at ( 0, 0 ) {};
\node[bla] at ( \aX, \aY) {};
\node[bll] at ( \aX\bX, \aY\bY) {};
\node[bll] (b11100) at ( \aX\bX\cX      , \aY\bY\cY      ) {};
\node[bla] (b11110) at ( \aX\bX\cX\dX   , \aY\bY\cY\dY   ) {};
\node[bla] (b01110) at (    \bX\cX\dX   ,    \bY\cY\dY   ) {};
\node[bll] at ( \cX\dX, \cY\dY) {};
\node[bll] at ( \dX, \dY) {};

\node[blo] at ( \bX, \bY) {};
\node[blo] at ( \cX, \cY) {};
\node[blo] at ( \aX\dX, \aY\dY) {};
\node[blo] at ( \bX\dX, \bY\dY) {};
\node[blo] at ( \aX\cX, \aY\cY) {};
\node[blo] at ( \bX\cX, \bY\cY) {};

 \node[blo] (b10110) at ( \aX   \cX\dX   , \aY   \cY\dY   ) {};
 \node[blo] (b11010) at ( \aX\bX   \dX   , \aY\bY   \dY   ) {};

 \node[vtx,label=left:$_{0000}$] (v00000) at ( 0              , 0               ) {\1}; 
 \node[vtx,label=below:$_{1000}$] (v00001) at ( \aX            , \aY             ) {\1}; 
 \node[vtx] (v00011) at ( \bX\aX         , \bY\aY          ) {\2}; 
 \node[vtx,label=above:$_{1110}$] (v00111) at ( \cX\bX\aX      , \cY\bY\aY       ) {\2}; 
 \node[vtx] (v01111) at ( \dX\cX\bX\aX   , \dY\cY\bY\aY    ) {\1}; 
 \node[vtx] (v01110) at ( \dX\cX\bX      , \dY\cY\bY       ) {\1}; 
 \node[vtx] (v01100) at ( \dX\cX         , \dY\cY          ) {\2}; 
 \node[vtx,label=below:$_{0001}$] (v01000) at ( \dX            , \dY             ) {\2}; 

 \node[vtx] (v00010) at ( \bX            , \bY             ) {}; 
 \node[vtx,label=below:$_{0010}$] (v00100) at ( \cX            , \cY             ) {}; 
 \node[vtx] (v00101) at ( \cX\aX         , \cY\aY          ) {}; 
 \node[vtx] (v00110) at ( \cX\bX         , \cY\bY          ) {}; 
 \node[vtx] (v01001) at ( \dX\aX         , \dY\aY          ) {}; 
 \node[vtx] (v01010) at ( \dX\bX         , \dY\bY          ) {}; 
 \node[vtx] (v01011) at ( \dX\bX\aX      , \dY\bY\aY       ) {}; 
 \node[vtx] (v01101) at ( \dX\cX\aX      , \dY\cY\aY       ) {}; 
 \draw[edge] (v00000) -- (v00010) -- (v00011) -- (v00001) -- (v00000);
 \draw[edge] (v00100) -- (v00110) -- (v00111) -- (v00101) -- (v00100);
 \draw[edge] (v01000) -- (v01010) -- (v01011) -- (v01001) -- (v01000);
 \draw[edge] (v01100) -- (v01110) -- (v01111) -- (v01101) -- (v01100);
 \draw[edge] (v00000) -- (v01000) -- (v01100) -- (v00100) -- (v00000);
 \draw[edge] (v00001) -- (v01001) -- (v01101) -- (v00101) -- (v00001);
 \draw[edge] (v00010) -- (v01010) -- (v01110) -- (v00110) -- (v00010);
 \draw[edge] (v00011) -- (v01011) -- (v01111) -- (v00111) -- (v00011);

 \end{tikzpicture}
\quad
 \raisebox{2em}{
 \begin{tikzpicture}[yscale=2.5,xscale=2.0]
\node[bll] at ( 0, 0 ) {};
\node[bla] at ( \aX, \aY) {};
\node[bla] at ( \aX\bX, \aY\bY) {};
\node[bll] (b11100) at ( \aX\bX\cX      , \aY\bY\cY      ) {};
\node[bll] (b11110) at ( \aX\bX\cX\dX   , \aY\bY\cY\dY   ) {};
\node[bla] (b01110) at (    \bX\cX\dX   ,    \bY\cY\dY   ) {};
\node[bla] at ( \cX\dX, \cY\dY) {};
\node[bll] at ( \dX, \dY) {};

\node[blo] at ( \bX, \bY) {};
\node[blo] at ( \cX, \cY) {};
\node[blo] at ( \aX\dX, \aY\dY) {};
\node[blo] at ( \bX\dX, \bY\dY) {};
\node[blo] at ( \aX\cX, \aY\cY) {};
\node[blo] at ( \bX\cX, \bY\cY) {};

 \node[blo] (b10110) at ( \aX   \cX\dX   , \aY   \cY\dY   ) {};
 \node[blo] (b11010) at ( \aX\bX   \dX   , \aY\bY   \dY   ) {};

 \node[vtx] (v00000) at ( 0              , 0               ) {\2}; 
 \node[vtx] (v00001) at ( \aX            , \aY             ) {\1}; 
 \node[vtx] (v00011) at ( \bX\aX         , \bY\aY          ) {\1}; 
 \node[vtx] (v00111) at ( \cX\bX\aX      , \cY\bY\aY       ) {\2}; 
 \node[vtx] (v01111) at ( \dX\cX\bX\aX   , \dY\cY\bY\aY    ) {\2}; 
 \node[vtx] (v01110) at ( \dX\cX\bX      , \dY\cY\bY       ) {\1}; 
 \node[vtx] (v01100) at ( \dX\cX         , \dY\cY          ) {\1}; 
 \node[vtx] (v01000) at ( \dX            , \dY             ) {\2}; 

 \node[vtx] (v00010) at ( \bX            , \bY             ) {}; 
 \node[vtx] (v00100) at ( \cX            , \cY             ) {}; 
 \node[vtx] (v00101) at ( \cX\aX         , \cY\aY          ) {}; 
 \node[vtx] (v00110) at ( \cX\bX         , \cY\bY          ) {}; 
 \node[vtx] (v01001) at ( \dX\aX         , \dY\aY          ) {}; 
 \node[vtx] (v01010) at ( \dX\bX         , \dY\bY          ) {}; 
 \node[vtx] (v01011) at ( \dX\bX\aX      , \dY\bY\aY       ) {}; 
 \node[vtx] (v01101) at ( \dX\cX\aX      , \dY\cY\aY       ) {}; 
 \draw[edge] (v00000) -- (v00010) -- (v00011) -- (v00001) -- (v00000);
 \draw[edge] (v00100) -- (v00110) -- (v00111) -- (v00101) -- (v00100);
 \draw[edge] (v01000) -- (v01010) -- (v01011) -- (v01001) -- (v01000);
 \draw[edge] (v01100) -- (v01110) -- (v01111) -- (v01101) -- (v01100);
 \draw[edge] (v00000) -- (v01000) -- (v01100) -- (v00100) -- (v00000);
 \draw[edge] (v00001) -- (v01001) -- (v01101) -- (v00101) -- (v00001);
 \draw[edge] (v00010) -- (v01010) -- (v01110) -- (v00110) -- (v00010);
 \draw[edge] (v00011) -- (v01011) -- (v01111) -- (v00111) -- (v00011);

 \end{tikzpicture}
 }
 \quad
 \begin{tikzpicture}[yscale=2.5,xscale=2.0]
\node[bla] at ( 0, 0 ) {};
\node[bll] at ( \aX, \aY) {};
\node[bll] at ( \aX\bX, \aY\bY) {};
\node[bla] (b11100) at ( \aX\bX\cX      , \aY\bY\cY      ) {};
\node[bla] (b11110) at ( \aX\bX\cX\dX   , \aY\bY\cY\dY   ) {};
\node[bll] (b01110) at (    \bX\cX\dX   ,    \bY\cY\dY   ) {};
\node[bll] at ( \cX\dX, \cY\dY) {};
\node[bla] at ( \dX, \dY) {};

\node[blo] at ( \bX, \bY) {};
\node[blo] at ( \cX, \cY) {};
\node[blo] at ( \aX\dX, \aY\dY) {};
\node[blo] at ( \bX\dX, \bY\dY) {};
\node[blo] at ( \aX\cX, \aY\cY) {};
\node[blo] at ( \bX\cX, \bY\cY) {};

 \node[blo] (b10110) at ( \aX   \cX\dX   , \aY   \cY\dY   ) {};
 \node[blo] (b11010) at ( \aX\bX   \dX   , \aY\bY   \dY   ) {};

 \node[vtx] (v00000) at ( 0              , 0               ) {\1}; 
 \node[vtx] (v00001) at ( \aX            , \aY             ) {\2}; 
 \node[vtx] (v00011) at ( \bX\aX         , \bY\aY          ) {\2}; 
 \node[vtx] (v00111) at ( \cX\bX\aX      , \cY\bY\aY       ) {\1}; 
 \node[vtx] (v01111) at ( \dX\cX\bX\aX   , \dY\cY\bY\aY    ) {\1}; 
 \node[vtx] (v01110) at ( \dX\cX\bX      , \dY\cY\bY       ) {\2}; 
 \node[vtx] (v01100) at ( \dX\cX         , \dY\cY          ) {\2}; 
 \node[vtx] (v01000) at ( \dX            , \dY             ) {\1}; 

 \node[vtx] (v00010) at ( \bX            , \bY             ) {}; 
 \node[vtx] (v00100) at ( \cX            , \cY             ) {}; 
 \node[vtx] (v00101) at ( \cX\aX         , \cY\aY          ) {}; 
 \node[vtx] (v00110) at ( \cX\bX         , \cY\bY          ) {}; 
 \node[vtx] (v01001) at ( \dX\aX         , \dY\aY          ) {}; 
 \node[vtx] (v01010) at ( \dX\bX         , \dY\bY          ) {}; 
 \node[vtx] (v01011) at ( \dX\bX\aX      , \dY\bY\aY       ) {}; 
 \node[vtx] (v01101) at ( \dX\cX\aX      , \dY\cY\aY       ) {}; 
 \draw[edge] (v00000) -- (v00010) -- (v00011) -- (v00001) -- (v00000);
 \draw[edge] (v00100) -- (v00110) -- (v00111) -- (v00101) -- (v00100);
 \draw[edge] (v01000) -- (v01010) -- (v01011) -- (v01001) -- (v01000);
 \draw[edge] (v01100) -- (v01110) -- (v01111) -- (v01101) -- (v01100);
 \draw[edge] (v00000) -- (v01000) -- (v01100) -- (v00100) -- (v00000);
 \draw[edge] (v00001) -- (v01001) -- (v01101) -- (v00101) -- (v00001);
 \draw[edge] (v00010) -- (v01010) -- (v01110) -- (v00110) -- (v00010);
 \draw[edge] (v00011) -- (v01011) -- (v01111) -- (v00111) -- (v00011);
 \end{tikzpicture}
 \quad
\raisebox{2em}{
 \begin{tikzpicture}[yscale=2.5,xscale=2.0]
\node[bll] at ( 0, 0 ) {};
\node[bll] at ( \aX, \aY) {};
\node[bla] at ( \aX\bX, \aY\bY) {};
\node[bla] (b11100) at ( \aX\bX\cX      , \aY\bY\cY      ) {};
\node[bll] (b11110) at ( \aX\bX\cX\dX   , \aY\bY\cY\dY   ) {};
\node[bll] (b01110) at (    \bX\cX\dX   ,    \bY\cY\dY   ) {};
\node[bla] at ( \cX\dX, \cY\dY) {};
\node[bla] at ( \dX, \dY) {};

\node[blo] at ( \bX, \bY) {};
\node[blo] at ( \cX, \cY) {};
\node[blo] at ( \aX\dX, \aY\dY) {};
\node[blo] at ( \bX\dX, \bY\dY) {};
\node[blo] at ( \aX\cX, \aY\cY) {};
\node[blo] at ( \bX\cX, \bY\cY) {};
 
 \node[blo] (b10110) at ( \aX   \cX\dX   , \aY   \cY\dY   ) {};
 \node[blo] (b11010) at ( \aX\bX   \dX   , \aY\bY   \dY   ) {};
 
 \node[vtx] (v00000) at ( 0              , 0               ) {\2}; 
 \node[vtx] (v00001) at ( \aX            , \aY             ) {\2}; 
 \node[vtx] (v00011) at ( \bX\aX         , \bY\aY          ) {\1}; 
 \node[vtx] (v00111) at ( \cX\bX\aX      , \cY\bY\aY       ) {\1}; 
 \node[vtx] (v01111) at ( \dX\cX\bX\aX   , \dY\cY\bY\aY    ) {\2}; 
 \node[vtx] (v01110) at ( \dX\cX\bX      , \dY\cY\bY       ) {\2}; 
 \node[vtx] (v01100) at ( \dX\cX         , \dY\cY          ) {\1}; 
 \node[vtx] (v01000) at ( \dX            , \dY             ) {\1}; 
 
 \node[vtx] (v00010) at ( \bX            , \bY             ) {}; 
 \node[vtx] (v00100) at ( \cX            , \cY             ) {}; 
 \node[vtx] (v00101) at ( \cX\aX         , \cY\aY          ) {}; 
 \node[vtx] (v00110) at ( \cX\bX         , \cY\bY          ) {}; 
 \node[vtx] (v01001) at ( \dX\aX         , \dY\aY          ) {}; 
 \node[vtx] (v01010) at ( \dX\bX         , \dY\bY          ) {}; 
 \node[vtx] (v01011) at ( \dX\bX\aX      , \dY\bY\aY       ) {}; 
 \node[vtx] (v01101) at ( \dX\cX\aX      , \dY\cY\aY       ) {}; 
 \draw[edge] (v00000) -- (v00010) -- (v00011) -- (v00001) -- (v00000);
 \draw[edge] (v00100) -- (v00110) -- (v00111) -- (v00101) -- (v00100);
 \draw[edge] (v01000) -- (v01010) -- (v01011) -- (v01001) -- (v01000);
 \draw[edge] (v01100) -- (v01110) -- (v01111) -- (v01101) -- (v01100);
 \draw[edge] (v00000) -- (v01000) -- (v01100) -- (v00100) -- (v00000);
 \draw[edge] (v00001) -- (v01001) -- (v01101) -- (v00101) -- (v00001);
 \draw[edge] (v00010) -- (v01010) -- (v01110) -- (v00110) -- (v00010);
 \draw[edge] (v00011) -- (v01011) -- (v01111) -- (v00111) -- (v00011);
 \end{tikzpicture}
 }
 }
 \caption{Four perfect $3$-colorings $f_{0,0}$, $f_{0,1}$, $f_{1,0}$, $f_{1,1}$,  of $Q_4$ with quotient matrix $[[1,1,2],[1,1,2],[1,1,2]]$}
 \label{f:H52} 
 \end{figure}%
} 

\title{On degree-$3$ and $(n-4)$-correlation-immune perfect colorings of $n$-cubes}
\author{Denis S. Krotov, Alexandr A. Valyuzhenich%
\thanks{D.S.K. and A.A.V. are with the Sobolev Institute of mathematics,
Novosibirsk 630090, Russia. E-mail: dk@ieee.org, graphkiper@mail.ru\\
This is the accepted version of the article D.\,S.\,Krotov, A.\,A.\,Valyuzhenich. On \mbox{degree-$3$} and $(n-4)$-correlation-immune
perfect colorings of $n$-cubes. 
\emph{Discrete Mathematics} 347(10):114138/1--14, 2024, doi
\href{https://doi.org/10.1016/j.disc.2024.114138}{10.1016/j.disc.2024.114138}
}}
\date{}

\begin{document}

\maketitle
\begin{abstract}
A perfect $k$-coloring of the Boolean hypercube $Q_n$ is a function from the set of binary words of length $n$ onto a $k$-set of colors such that for any colors $i$ and $j$ every word of color $i$ has exactly $S(i,j)$ neighbors (at Hamming distance $1$) of color $j$, where the coefficient $S(i,j)$ depends only on $i$ and $j$ but not on the particular choice of the word. The $k$-by-$k$ table of all coefficients $S(i,j)$ is called the quotient matrix. We characterize perfect colorings of $Q_n$ of degree at most~$3$, that is, with quotient matrix whose all eigenvalues are not less than $n-6$, or, equivalently, such that every color corresponds to a Boolean function represented by a polynomial of degree at most~$3$ over~$R$. Additionally, we characterize $(n-4)$-correlation-immune perfect colorings of~$Q_n$, whose all colors correspond to $(n-4)$-correlation-immune Boolean functions, or, equivalently, all non-main (different from $n$) eigenvalues of the quotient matrix are not greater than $6-n$.

\end{abstract}

\section{Introduction}




An arbitrary surjective function 
from the vertex set
of a graph~$\Gamma$ onto a finite set~$K$ 
(of \emph{colors})
of cardinality~$k$ is called 
a (vertex) \emph{coloring}, 
or \emph{$k$-coloring}, of~$\Gamma$. 
A $k$-coloring of a graph is called \emph{perfect}
if there is a $k$-by-$k$ 
matrix $\{S_{i,j}\}_{i,j\in K}$
(the \emph{quotient matrix})
such that every vertex of color~$i$
has exactly~$S_{i,j}$ neighbors
of color~$j$.
An {\em eigenvalue} of a perfect coloring is an eigenvalue of its quotient matrix.
It is well known  that if $f$ is a perfect coloring of a graph~$\Gamma$, then every eigenvalue of $f$ is an eigenvalue of the adjacency matrix of~$\Gamma$.
Two colorings $f$ and $g$ of a graph~$\Gamma$ 
are \emph{equivalent} if $g(\cdot)\equiv\pi(f(\alpha(\cdot)))$, 
where $\alpha$ is an automorphism of~$\Gamma$ 
and $\pi$ is a bijection between the color sets of~$f$ and~$g$.

The {\em $n$-dimensional hypercube} $Q_n$ is defined as follows.
The vertex set of~$Q_n$ is $\mathbb{Z}_{2}^n$, and two vertices are adjacent if they differ in exactly one coordinate.
This graph has $n+1$ distinct eigenvalues $\lambda_i(n)=n-2i$, $0\leq i\leq n$, with the corresponding eigenspaces~$U_i(n)$.
Denote by $U(n)$ the set of all real-valued functions defined on the vertex set of~$Q_n$.
Every function $f\in U(n)$ is uniquely represented as the sum
\begin{equation}\label{eq:sm}
    f=\sum_{i=0}^{n}f_i,\qquad \mbox{where } f_i\in U_i(n).
\end{equation}
A function $f\in U(n)$ has {\em degree}~$d$ if $f_d\not\equiv 0$ and $f_i\equiv 0$ for all $i>d$ in~\eqref{eq:sm}.
A perfect coloring of $Q_n$ has {\em degree}~$d$ if its smallest eigenvalue is $n-2d$.
A function $f$ on $\mathbb{Z}_{2}^n$ is said to depend \emph{essentially} in the $i$th argument, $1\leq i\leq n$, if there are two vectors $\vc x,\vc y\in \mathbb{Z}_{2}^n$ differing only in the $i$th position such that $f(\vc x)\neq f(\vc y)$.

Perfect colorings of distance-regular graphs have been extensively investigated by various researchers under different names.
In particular, perfect $2$-colorings of small degree were studied for Hamming graphs \cite{Meyer:cyccond,MV:eq2eigen}, Johnson graphs \cite{EGGV:Johnson,GavGor:2013,Vorobev:equit_Johnson:arx},  Grassmann graphs \cite{GavMat:2018,GavMog:2014,Matkin:2018}.
Boolean functions of degree~$1$ on several classical association schemes were investigated in \cite{FilmusIhringer:2019}, and
Boolean functions of degree~$2$ on Grassmann graphs were investigated in~\cite{dBDIM:2023}.
Now, we briefly discuss results on perfect $2$-colorings of~$Q_n$ of degree at most~$3$.
The famous Nisan--Szegedy theorem~\cite{NisSze:94} states that any Boolean function of degree~$d$ has at most $d\cdot 2^{d-1}$ essential variables
(this bound was improved to $6.614 \cdot 2^d$ and then to $4.416\cdot 2^d$ in~\cite{CHS:2020} and~\cite{Wellens:tighter}, respectively).
This immediately implies that any perfect $2$-coloring of~$Q_n$ of degree~$d$ has at most $d\cdot 2^{d-1}$ essential variables.
So, we have the following for $d\leq 2$:

\begin{itemize}
    \item Any perfect $2$-coloring of $Q_n$ of degree $1$ has at most one essential variable. 
    Therefore, there is exactly one such perfect $2$-coloring up to equivalence.
    Note that this fact is a very special case of the result obtained by Meyerowitz in \cite{Meyer:cyccond}
    (Meyerowitz characterized all completely regular codes in the Hamming graph containing the second largest eigenvalue of the graph in their spectrum).

\item Any perfect $2$-coloring of~$Q_n$ of degree~$2$ 
has at most~$4$ essential variables. 
So, all such colorings can be found by enumeration of all $2$-colorings of~$Q_4$.
In particular, there are three distinct degree-$2$ perfect $2$-coloring of~$Q_n$ up to equivalence {(see Fig.~\ref{f:F0})}. 
\figureA
In other terms, such colorings were characterized by Camion et al. in~\cite{CCCS:92}. 
In addition, Mogilnykh and Valyuzhenich~\cite{MV:eq2eigen} showed that any 
degree-$2$
perfect 
$2$-coloring of the Hamming graph 
either has at most three essential variables
or is obtained by one of two special constructions.
\end{itemize}

Taking into account the known 
connection between low-degree
Boolean functions and 
high-order correlation-immune Boolean functions (see Lemma~\ref{l:deg-ci} below), the classification of perfect $2$-colorings of $Q_n$ of degree at most~$3$
can be derived from the classification of balanced $(n-4)$-correlation-immune functions recently obtained by Rasoolzadeh~\cite{Rasoolzadeh2023}.
In this work, we classify perfect colorings of $Q_n$ of degree at most~$3$ with $2$ or more colors.
Moreover, motivated by the above-mentioned connection with correlation-immune functions, 
we classify also
the perfect colorings
whose all non-main (different from~$n$)
eigenvalues are not larger than $6-n$,
which means that these colorings are 
$(n-4)$-correlation-immune.

The paper is organized as follows.
Section~\ref{s:bg} contains the background material, including the basic facts we use and necessary definitions.
In Section \ref{s:k:n-6}, we classify
perfect colorings of $Q_n$ of degree at most $3$.
We describe an algorithm
that derives
all such colorings of $Q_{10}$ 
from the known characterization of $(n-4)$-resilient functions~\cite{Rasoolzadeh2023},
show that such colorings of $Q_n$, $n>10$, have nonessential arguments,
and describe the structure
of all such colorings with more
than~$2$ colors.
In Section \ref{s:6-n}, we classify
the $(n-4)$-correlation-immune perfect colorings of~$Q_n$.
\section{Background}\label{s:bg}


\subsection{Density, correlation-immune 
and resilient functions}

For a coloring of a finite graph, 
we define the 
\emph{density}~$\rho(i)$
of each color~$i$ as the proportion
of the vertices of color~$i$
among the all vertices of the graph.
The \emph{density vector} of 
a coloring is the list of densities 
of all colors;
it has sum~$1$.
For a perfect coloring,
the density vector is uniquely determined
by the quotient matrix~$\{S_{i,j}\}_{i,j}$,
because 
\begin{equation}\label{eq:rhoij}
\rho_i S_{i,j} = \rho_j S_{j,i}
\end{equation}
by double-counting the edges between 
color-$i$ and color-$j$ vertices. 

A coloring of $Q_n$
(in particular, a Boolean function,
as a $1$- or $2$-coloring of~$Q_n$)
is called 
{$t$-correlation-immune}
if the density vector
of its retract coloring 
on a subgraph isomorphic to 
$Q_{n-t}$ does not depend 
on the choice of such subgraph
(note that every such subgraph 
is induced by the vertices 
with fixed values in some fixed $t$ positions).
A $k$-coloring of $Q_n$
(in particular, a Boolean function if $k=2$) is called \emph{$t$-resilient}
if the density vector
of its retract coloring 
on every subgraph isomorphic to 
$Q_{n-t}$ is $(\frac1k,\ldots,\frac1k)$.

 The following lemma contains 
 three known important observations
 connecting Boolean functions
 of high-order resilience and of low degree.
 The first one can be found, e.g., in
 \cite{Donnell:2014} (the remark between Propositions~6.23 and~6.24); the last two are straightforward. 
 The $i$th argument of a Boolean function~$f$ is called \emph{linear} for~$f$ if $f(x_1,\ldots,x_{i-1},x_{i}+1,x_{i+1},\ldots,x_n)\equiv f(x_1,\ldots,x_n)+1$.
\begin{lemma}\label{l:deg-ci} 
Assume that $f$ and $f'$ are Boolean functions
on $Q_n$ whose values coincide on one bipartite part of 
$Q_n$ and are opposite on each vertex of the other bipartite part. 
\begin{itemize}
    \item[\rm(a)]
The degree of $f$ does not exceed~$d$ if and only if 
$f'$ is $(n-d-1)$-resilient.  
    \item[\rm(b)]
Moreover, $f$ is a perfect coloring 
with quotient matrix 
$\begin{pmatrix}a&b\\b&a
\end{pmatrix}$
if and only if 
$f'$ is a perfect coloring 
with quotient matrix
$\begin{pmatrix}b&a\\a&b
\end{pmatrix}$,
$a,b>0$.
    \item[\rm(c)]
    The $i$th argument is linear for~$f'$
    if and only if $f$ does not depend on its value.
\end{itemize}
\end{lemma}

We note that claim~(a)
of Lemma~\ref{l:deg-ci} 
is valid for $d=n$
as well because 
by the definition 
every Boolean function is
$-1$-resilient.
\begin{lemma}[folklore, see e.g. {\cite[Proposition~1]{Kro:OA13}}]\label{l:perf-ci} 
A non-constant perfect coloring of $Q_n$ 
is $t$-correlation-immune,
where $\theta = n-2(t+1)$ 
is the second largest eigenvalue of the quotient matrix,
$t = \frac{n-\theta}2-1$.
\end{lemma}

Fon-Der-Flaass \cite{FDF:CorrImmBound} established the following bound on the correlation immunity of non-constant unbalanced Boolean functions.

\begin{lemma}\label{L:CorrImmBound}
For a $t$-correlation-immune Boolean
function in $n$ arguments with density of~$1$ different
from $0$, $\frac 12$, $1$, it holds
\begin{equation}
    \label{eq:c-i}
t \le \frac{2n}{3}-1.
\end{equation}
In the case of equality,
the function is a perfect $2$-coloring with eigenvalue~$-\frac{n}{3}$.    
\end{lemma}
\begin{corollary}\label{c:FDF}
Any $t$-correlation-immune $k$-coloring of~$Q_n$, $k\ge 3$,
satisfies~\eqref{eq:c-i}.
In the case of equality,
the coloring is perfect, 
with exactly two eigenvalues~$n$
and~$-\frac{n}{3}$.
\end{corollary}
\begin{proof}
Since $k\ge 3$, at least $k-1$ colors
have density strictly
between $0$ and~$\frac{1}{2}$.
Replacing one of these colors by~$1$
and the other colors by~$0$, we obtain
a Boolean function satisfying 
the hypothesis
in Lemma~\ref{L:CorrImmBound}.
Hence,
\eqref{eq:c-i} holds and, moreover,
in the case of equality each of these $k-1$ functions is a 
perfect $2$-coloring with eigenvalue~$-\frac{n}3$ 
and by Corollary~\ref{l:union}
the original $k$-coloring
is perfect with two 
eigenvalues $n$ and~$-\frac{n}3$.
\end{proof}
\begin{corollary}\label{c:FDF2}
    If there exists a perfect coloring of $Q_n$ with exactly
    two eigenvalues $n$ and~$\theta$,
    then either it is a perfect $2$-coloring with symmetric quotient matrix, or $\theta \geq -\frac{n}{3}$.    
\end{corollary}
\begin{proof}
It follows from Lemmas \ref{l:perf-ci}, \ref{L:CorrImmBound}, and Corollary \ref{c:FDF}.  
\end{proof}

In the rest of this subsection,
we describe the known results on 
the characterization of 
high-order 
resilient Boolean functions, which play a key role 
in our current study.

In~\cite{TarKir2000}, 
Tarannikov and Kirienko 
showed that for every 
positive integer~$m$, 
there exists a minimum number~$p(m)$
such that for $n > p(m)$,
any $(n-m)$-resilient
$n$-argument function
has at least $n-p(m)$
linear arguments; 
in particular, the total number
of such functions is a polynomial
of degree~$p(m)$ in~$n$.
They also 
proved that $p(4)=10$
(a computer-free prove was later suggested by Zverev~\cite{Zverev:2008})
and 
developed a technique,
based on the Walsh--Hadamard
transform,
to analyse such functions.
Utilizing that technique,
Kirienko~\cite{Kirienko:2004} 
derived the exact number 
\begin{multline}\label{eq:Kir}
  \frac{1}{2} n^{10 }
+ \frac{7}{6} n^9 
+ \frac{890}{9} n^8 
- \frac{10903}{9} n^7 
+ \frac{64288}{45} n^6 \\
+ \frac{953308}{45} n^5
- \frac{1341569}{18} n^4 
+ \frac{899251}{18} n^3 
+ \frac{365018}{5} n^2 
- \frac{1048961}{15} n + 2  
\end{multline} 
of $(n-4)$-resilient functions 
in $n$ arguments, for any~$n$.
Recently,
Rasoolzadeh~\cite{Rasoolzadeh2023}
characterized all such functions
up to equivalence
and provided a representative 
for each equivalence class.
By counting the sizes
of the equivalence classes found in~\cite{Rasoolzadeh2023},
one can check that the two results obtained independently with different techniques
are in agreement with each other (it is sufficient to compare the values of~\eqref{eq:Kir} for $n=0,\ldots,10$:
2, 4, 16, 256,
     12870,
     807980,
    16750860,
    126113920,
   605047818, 
   2220784820,
  6799438888).

For the completeness 
of this mini-survey,
we mention the earlier results
with the characterization 
of $(n-2)$- and $(n-3)$-resilient functions, among other results in~\cite{Siegenthaler:84} and~\cite{CCCS:92}, respectively;
and the characterization of 
orthogonal arrays for many series
of parameters~\cite{SEN:2010:OA}, 
which includes the parameters of 
all unbalanced 
$(n-4)$-correlation-immune
Boolean functions (by Lemma~\ref{L:CorrImmBound}, $n\le 9$; the case $n=9$ was essentially characterized in~\cite{Kirienko2002}).
So, we can say that~\cite{Rasoolzadeh2023}
completes the characterization
of the $(n-4)$-correlation-immune
Boolean functions by solving the balanced case.

\begin{remark}\label{r:deg3}
   By Lemma~\ref{l:deg-ci},
expression~\eqref{eq:Kir}
gives also the value
of Boolean functions
of degree at most~$3$.
Moreover, representatives
of the equivalence classes
of such functions can be
obtained from 
the representatives
of $(n-4)$-resilient functions
presented 
in~\cite{Rasoolzadeh2023} 
by inverting the values
of the function
in one of the bipartite parts
of~$Q_n$.
\end{remark}

\subsection{Properties of perfect colorings, the coarsest equitable refinement}

A $\{0,1\}$-valued function~$t$ on
the vertex set of a graph~$\Gamma$
is called a \emph{fiber} of 
a vertex coloring~$f$ of~$\Gamma$ if
$t^{-1}(1)=f^{-1}(i)$ 
for some color~$i$, 
i.e., $t$ is the characteristic
function of the $f$-preimage of~$i$.

\begin{lemma}[see e.g.~{\cite[Ch.~5, Lemma~2.2]{Godsil93}}]\label{l:eig}
For a perfect coloring~$f$ of a graph~$\Gamma$, 
every fiber
is the sum of eigenfunctions 
of~$\Gamma$
corresponding to eigenvalues 
of the quotient matrix of~$f$.
Moreover, every eigenvalue 
of the quotient matrix 
is an eigenvalue of~$\Gamma$.
\end{lemma}

\begin{corollary}
\label{l:union}
\begin{itemize}
    \item[\rm (a)]  If a perfect $k$-coloring $f$ of an $r$-regular graph
has only two eigenvalues $r$ and~$\theta$,
then merging any two colors
results in a perfect $(k-1)$-coloring.
   \item[\rm (b)]
Inversely, if there are a perfect $k$-coloring~$f$
and a perfect $2$-coloring~$g$
of the same regular graph, their quotient 
matrices have only two eigenvalues~$r$ and~$\theta$,
and the coloring 
$h(\vc{x})=(f(\vc{x}),g(\vc{x}))$
is a $(k+1)$-coloring,
then it is a perfect
coloring with two eigenvalues~$r$ and~$\theta$.
\end{itemize}
\end{corollary}
\begin{proof}
(a)
By Lemma~\ref{l:eig},
each (say, the $i$th) fiber of $f$ is the sum of 
a constant (an eigenfunction with eigenvalue~$r$) and an eigenfunction with eigenvalue~$\theta$. 
It is straightforward to see that 
such a fiber is a perfect $2$-coloring with quotient matrix
$\Quot{a_i& b_i\\c_i&d_i}$
where $a_i+b_i=c_i+d_i=r$
and $a_i-c_i=d_i-b_i=\theta$,
i.e.,
$\Quot{r-b_i&b_i \\ r-b_i-\theta & b_i+\theta}$
for some~$b_i$. It is now straightforward from the definition that $f$ has quotient matrix $(S_{i,j})$
where $S_{i,j} = b_j$ if
$i\ne j$. Claim~(a) can now
be directly checked.

(b) Obviously, each fiber $f_1$, \ldots, $f_k$ of~$f$ and each fiber $g_1$, $g_2$ of~$g$ is a linear combination
of fibers of~$h$. On the other hand, the fibers $f_1$, \ldots, $f_k$, $g_1$ are linearly independent.
It follows that each of $k+1$
fibers of~$h$ is a linear combination of fibers of~$f$ and~$g$ and hence is a perfect $2$-coloring.
Similarly to~(a), $h$ is a perfect $(k+1)$-coloring with quotient matrix described as in~(a).
\end{proof}

For two vertex colorings $f$ and~$g$
of the same graph, it is said
that $g$ is 
a \emph{refinement} of $f$,
or $f$ is \emph{coarser} than~$g$, or $g\preceq f$, 
if $f(\vc x)\ne f(\vc y)$ implies $g(\vc x)\ne g(\vc y)$ (equivalently, every fiber of~$f$ is the sum of some fibers of~$g$).

It is well known, 
see e.g.~\cite{McKay:80}, 
that for every vertex coloring~$f$
there is a 
\emph{coarsest equitable refinement}
of~$f$,
i.e., a perfect coloring~$g$
of the same graph such that 
$g \preceq f$ and 
for any perfect coloring~$h$,
$h \preceq f$ implies $h \preceq g$.

We will use the following two simple known facts (the first one is trivial and needs no proof).
\begin{lemma}\label{l:fiber}
Every fiber $t$ of a perfect coloring
is a fiber of
the coarsest equitable refinement of~$t$.
\end{lemma}
\begin{lemma}
\label{l:spe<}
If a perfect coloring $g$ is a refinement of a perfect coloring~$f$,
then all eigenvalues of~$f$ are eigenvalues of~$g$.
\end{lemma}
\begin{proof}[A sketch of proof]
 To avoid complicate indices, 
 consider an example where $f$ and~$g$ are perfect $3$-
 and $5$-colorings with quotient matrices $S=\{S_{i,j}\}_{i,j\in\{a,b,c\}}$ 
 and $T=\{T_{i,j}\}_{i,j\in\{a',b',b'',c',c''\}}$, respectively, and
 \begin{equation}\label{eq:1,23,45}
   f^{-1}(a)=g^{-1}(a'), \quad
 f^{-1}(b)=g^{-1}(b')\cup g^{-1}(b''), \quad
 f^{-1}(c)=g^{-1}(c')\cup g^{-1}(c'').
 \end{equation}
 From~\eqref{eq:1,23,45} and the definition of perfect coloring,
 we see the following:
 $$
 \newcommand{\Sss}{\phantom{S_{a,a}}}
 \begin{array}{l|l|l}
 S_{a,a}=T_{a',a'},   & S_{a,b}=T_{a,b'}+T_{a,b''}, &  S_{a,c}=T_{a,c'}+T_{a,c''},   \\ \hline
 S_{b,a}=T_{b',a'}  & S_{b,b}=T_{b',b'}+T_{b',b''}   &  S_{b,c}=T_{b',c'}+T_{b',c''} \\
  \Sss=T_{b'',a'}, &   \Sss=T_{b'',b'}+T_{b'',b''}, &   \Sss=T_{b'',c'}+T_{b'',c''},\\  \hline
 S_{c,a}=T_{c',a'}  & S_{c,b}=T_{c',b'}+T_{c',b''}  &  S_{c,c}=T_{c',c'}+T_{c',c''} \\
  \Sss=T_{c'',a'}, &  \Sss=T_{c'',b'}+T_{c'',b''}, &   \Sss=T_{c'',c'}+T_{c'',c''}.
 \end{array}
$$
 It is now straightforward to see that an eigenvector
 $(x_a,x_b,x_c)$ of~$S$ corresponds to the eigenvector
 $(x_{a},x_{b},x_{b},x_{c},x_{c})$ of~$T$ with the same eigenvalue.
 In the general case, the proof is similar.
\end{proof}

Finally, the following
fact is also a known corollary
of the property of perfect colorings called
distance invariance, 
see e.g.~\cite{Kro:struct}: if we know the color of a vertex~$\vc x$, 
then we can calculate
 the number of vertices 
 of each color at distance~$i$ from~$\vc x$, $i=1,2,\ldots,n$.
\begin{lemma}\label{l:antipod}
For a perfect coloring of~$Q_n$,
the color of a vertex~$\vc x$ 
is uniquely determined
by the color of the complement vertex $\vc x + \overline 1$.
In particular there is a matching on the set of colors
such that any two complement vertices have matched colors
(each color is matched with another color or with itself).
Two matched colors have equal densities.
\end{lemma}


\section
[Perfect colorings of Qn with 
smallest eigenvalue n-6]
{Perfect colorings 
of $Q_n$ of degree at most $3$}
\label{s:k:n-6}
In this section,
we consider the classification
of perfect colorings
of $Q_n$ of degree at most $3$.
For brevity, we will call them 
\emph{low-degree} perfect colorings.
We use the fact that
every fiber of such a coloring
is a Boolean function
of degree at most~$3$,
such a function cannot
have more than~$10$
essential arguments,
and all such functions
have been characterized
up to equivalence in
terms of $(n-4)$-resilient
functions
in~\cite{Rasoolzadeh2023},
see Remark~\ref{r:deg3}.
In Subsection~\ref{s:class},
we describe
the computer-aided 
classification of 
low-degree
perfect colorings of $Q_{10}$
and prove that every 
low-degree perfect coloring
of $Q_n$ with $n>10$
has a nonessential argument.
In Subsection~\ref{s:descr},
we describe constructions
that cover all found low-degree
perfect colorings with more than~$2$ colors.

\subsection{Classification}
\label{s:class}
Here, we describe the algorithm used
to find representatives of all equivalence classes
of perfect colorings of~$Q_n$ of degree at most~$3$. 
It
was implemented
for $n=10$, 
and the additional
check described
in
Proposition~\ref{p:n<11}
convinces us that this algorithm
does not give more perfect colorings
if it runs with any larger value 
of~$n$.
\begin{itemize}
    \item[(i)] From \cite{Rasoolzadeh2023},
    we know all nonequivalent
    $(n-4)$-resilient 
    Boolean functions in $n$ arguments, for any~$n$ (for $n\le 10$, the number of equivalence classes and representatives are 
    explicitly given in \cite{Rasoolzadeh2023}; 
    for $n > 10$, the number 
    of equivalence classes
    is the same as for $n=10$, i.e., $961$).
    Using Lemma~\ref{l:deg-ci},
    we find all nonequivalent 
    Boolean functions
    ($2$-colorings)
    of degree at most~$3$.
    In particular, 
    independently on~$n$,
    such functions have at most~$10$ essential arguments.
    
    \item[(i')]
    For each function from 
    the step above, we find 
    its coarsest equitable refinement
    and if this refinement 
    is low-degree, 
    collect all fibers of all such 
    refinements.
    For each such fiber,
    by applying all automorphisms
    of~$Q_n$, 
    we find the equivalence class and
    denote the union of all
    such equivalence classes by~$T$.
    By Lemmas~\ref{l:fiber} 
    and~\ref{l:spe<}, $T$ 
    is the set of all fibers
    of low-degree perfect colorings
    of~$Q_n$.

    \item[] We start with 
    the collection~$P$ of perfect colorings
    consisting of only one
    constant $1$-coloring
    and iterate steps~(ii) and~(iii) below while step~(iii) 
    discovers new equivalence classes
    of perfect colorings.
    \item[(ii)] 
    For each representative~$g$ of a perfect $k$-coloring from the previous step and for each fiber~$t$ from~$T$
    such that $g$ is constant
    on the set of ones of~$t$
    but $t$ is not a fiber of~$g$, we consider the 
    coloring
    $h(x) = (g(x),t(x))$
    (which necessarily has $k+1$ colors)
    and keep
    nonequivalent representatives 
    of such colorings in a temporary
    collection~$S$.
    \item[(iii)] For each colorings 
    from~$S$,
    we find  
    the coarsest equitable refinement
    and check its eigenvalues.
    We restock $P$ by
    representatives
    of new found equivalence classes 
    of perfect colorings with all 
    eigenvalues not less than~$n-6$. 
    \item[] If step~(iii)
    produces no new colorings,
    i.e., does not increase~$P$,
    then the algorithm stops with the output~$P$.
    Otherwise, go to~(ii).
\end{itemize}
We ran the algorithm with $n=10$. Steps~(i), (i'), and~(iii) were implemented in \texttt{sage}~\cite{sage},
which has a convenient library
of graph tools (in particular,
a tool 
to find the 
coarsest equitable refinement).
Step~(ii) was implemented in \texttt{c++},
because the speed was 
important to run
over the large amount 
of Boolean functions in~$T$.
All calculations took 
several hours of CPU time.
As an intermediate result, the set~$T$
contained $2$~fibers with $128$ ones (with $3$ and $4$ essential arguments),
$8$~fibers of with $256$ ones
(with, $2$, $3$, $3$,
$4$, $5$, $5$, $6$, $6$ 
essential arguments),
$2$~fibers of with $384$ ones 
($3$, $4$),
$36$~fibers of with $512$ ones 
($1$, $2$, $3$, $4$, $5$, $7\times6$, $12\times7$, $9\times8$, $2\times9$, $10$),
$1$~fiber of with $768$ ones 
($3$ essential arguments), 
and $1$~constant fiber with $1024$ ones.
The first iteration of steps~(ii)--(iii) 
found $37$ nonequivalent low-degree perfect $2$-colorings,
five $3$-colorings, and four $3$-colorings. 
The $2$nd iteration 
found $2$ new perfect $3$-colorings,
$44$ new $4$-colorings,
three $6$-colorings,
and five $8$-colorings.
The $3$rd iteration found
$3$ new $4$-colorings.
The $4$th did not find new colorings.

\begin{proposition}\label{p:alg}
    Every perfect coloring of $Q_n$ of degree at most~$3$
    is obtained, up to equivalence,
    by the algorithm described above.
\end{proposition}
\begin{proof}
Assume that
after some $i$th iteration
of the algorithm,
there is a low-degree perfect coloring $f$ 
of $Q_n$ 
that is not found up to equivalence
(at step~(iii) of the algorithm).
It is sufficient to show 
that 
\begin{itemize}
    \item[(*)]
\emph{the next iteration of the algorithm
will find at least one new equivalence class of low-degree perfect colorings.}
\end{itemize}
This means that when algorithm stops
(it cannot run forever because the number of perfect colorings is finite), all low-degree perfect colorings are found, which is exactly
what the proposition claims.

Let $g$ 
be a perfect coloring such that
\begin{itemize}
    \item[(a)] $g\succeq f'$ for some $f'$ equivalent to~$f$,
    \item[(b)] $g$ was found 
    during the first~$i$ iterations of the algorithm (i.e., belongs to~$P$ after the $i$th iteration),
    \item[(c)] $g$ has the maximum number~$s$ of colors among the colorings satisfying~(a) and~(b)
(the set of such colorings is not empty
because it contains
the constant $1$-coloring).
\end{itemize}
By our hypothesis on~$f$, we have
$f' \ne g$, and hence there is 
a fiber~$t$ of~$f'$
that is not a fiber of~$g$. 
It is of degree at most~$3$
because $f$ is low-degree;
therefore, 
it belongs to the set~$T$
found at
step~(i') of the algorithm.
Denote by~$h$ the coloring 
$h(x)=(g(x),t(x))$.
It is an $(s+1)$-coloring;
moreover, $f'\preceq h \prec g$.
Denote by~$h'$ the coarsest equitable refinement
of~$h$. We have
$f'\preceq h'\preceq h \prec g$,
and by Lemma~\ref{l:spe<}
$h'$ is low-degree.
We see that the perfect coloring~$h'$ is obtained
from~$g$ at the $(i+1)$th iteration of the algorithm 
(step~(ii) generates~$h$,
then step~(iii) finds~$h'$),
$h'$ is coarser than~$f'$
and has more colors than~$g$.
By our hypothesis on~$g$, 
the perfect coloring~$h'$ new,
i.e., was not found during the first $i$ iterations.
This completes the proof of~(*).
\end{proof}


Next, we show that 
the classification of
low-degree perfect colorings
in $Q_{10}$
actually completes
the classification of
low-degree perfect colorings
in~$Q_n$ for any~$n$.

\begin{proposition}
\label{p:n<11}
For $n>10$, every 
perfect coloring of~$Q_n$ of degree at most~$3$
has at most~$10$ 
essential arguments.
\end{proposition}
From the results for Boolean functions~\cite{TarKir2000},
we know that every color component
has at most~$10$ 
essential arguments, but this does not directly imply
that these essential arguments are common for all colors
of the coloring.
Unfortunately, our proof
of Proposition~\ref{p:n<11} also includes
 computing;
it would be interesting
to find a purely theoretical
explanation, maybe even for all colorings (not only perfect), generalizing the proof~\cite{Zverev:2008} of the two-color case.
\begin{proof}
Seeking a contradiction,
we assume that there is
a low-degree 
perfect coloring~$f$
in~$Q_n$
without nonessential
arguments,
where $n$ is larger than~$10$.
Let us consider theoretically
the run of the algorithm 
described above 
for perfect colorings 
of~$Q_n$.
By Proposition~\ref{p:alg},
we know that $f$ 
is obtained
at some iteration
of the algorithm.
On the other hand,
the starting point 
of the algorithm 
is the Boolean
functions of degree 
at most~$3$, which have
at most~$10$ essential arguments, 
see step~(i).
This means that the number
of essential arguments
has been increased, exceeding~$10$
at some step.
It cannot be step~(iii) because 
of the following 
simple observation:
\begin{itemize}
    \item[] Claim~(*). \emph{The coarsest equitable refinement
    of a vertex coloring~$h$ of~$Q_n$ 
    does not have more essential arguments than $h$ has.}
    \item[] \emph{Proof of~(*).} Assume without loss of generality that $h$ does not depend essentially on the 
    first argument, i.e., $h=\tau h$, where
    $\tau h(x_1,x_2,...,x_{n})=h(1-x_1,x_2,...,x_{n})$.
    If $f$ is the coarsest equitable refinement of~$h$,
    then, obviously, the coarsest equitable refinement of~$\tau h$ is~$\tau f$. 
    Since $h=\tau h$, we see that
    $f$ and $\tau f$ are the same.
    So, $f$ also nonessentially depends on the
    first argument.
\end{itemize}
So, the number of essential arguments can firstly exceed $10$ at step~(ii) only,
where two colorings~$g$
and~$t$ are combined
resulting 
in a new coloring~$h$, 
$h(x) = (g(x),t(x))$,
which can potentially
have more essential arguments
than each of~$g$ and~$t$.
We assume that each of~$g$ and~$t$ 
has at most~$10$
essential arguments,
while $h$ has~$11$ or more
essential arguments
(the last means that there is an argument essential for~$g$
and nonessential for~$t$ 
and an argument essential for~$t$
and nonessential for~$g$).
By the condition of step~(ii),
$g$ is constant 
on the set of ones of~$t$,
and we assume without
loss of generality
that this constant is~$0$,
i.e., the sets of nonzeros 
of~$g$ and~$t$ do not intersect.
The colorings~$g$ and~$t$
satisfy the hypothesis 
of the following auxiliary 
claim.
\begin{itemize}
    \item[] Claim~(**). \emph{Let
    $g$ and $t$ be
    two colorings of $Q_n$ with disjoint 
    sets of nonzeros.
    We assume that there are
    two numbers in $\{1,\ldots,n\}$, $i$ and~$j$, such that
    $g$ depends essentially 
    in the $i$th argument
    and nonessentially 
    in the $j$th argument,
    while $t$, inversely,
    depends nonessentially 
    in the $i$th argument
    and essentially 
    in the $j$th argument.
    Denote by $\sigma_{ij}$ the operator
    that swaps the $i$th and $j$th arguments of a function:
    \begin{multline*}
    \sigma_{ij} t (x_1,...,x_{i-1},x_{i},x_{i+1},...,x_{j-1},x_{j},x_{j+1},...,x_n)\\
    =t(x_1,...,x_{i-1},x_{j},x_{i+1},...,x_{j-1},x_{i},x_{j+1},...,x_n).    
    \end{multline*}
        We claim that $g$ and $\sigma_{ij}t$
    also have disjoint sets of
    nonzeros. 
    Moreover, the number of essential arguments
    of the coloring
    $h_{ij}(x)=(g(x),\sigma_{ij}t(x))$ 
    is one smaller than
    that number for
    $h(x)=(g(x),t(x))$. }
    \item[] \emph{Proof of~(**).} The second statement (about 
    the number of essential arguments) is obvious; let us prove the first one. 
    Denote by $\tau_j$
    the operator inverting the $j$th argument
    of a function:
    $$\tau_j t(x_1,...,x_{j-1},x_{j},x_{j+1},...,x_n)=t(x_1,...,x_{j-1},1-x_{j},x_{j+1},...,x_n).$$
    Since $g$ and~$t$ have no common nonzeros
    and $g=\tau_j g$,
    we see that $g$ and~$\tau_j t$
     have no common nonzeros as well.
    It follows that $g$ and~$t_{\mathrm{mx}}$
    have no common nonzeros, 
    where 
    $t_{\mathrm{mx}}(x)=\max\{t(x),\tau_j t(x)\}$.
    Readily, $t_{\mathrm{mx}}$  
    depends nonessentially
    on the both $i$th and 
    $j$th arguments, which yields
    $t_{\mathrm{mx}}= \sigma_{i,j} t_{\mathrm{mx}}$. 
    It follows that  
    $t_{\mathrm{mx}}$ majorizes $\sigma_{i,j} t$:
    $$t_{\mathrm{mx}}(x) = \sigma_{i,j} t_{\mathrm{mx}}(x) 
    = \max \{\sigma_{i,j} t(x),\sigma_{i,j} \tau_j t(x) \}
    \ge \sigma_{i,j} t(x)\quad\mbox{for all~$x$.}$$ 
    Therefore, $g$ and $\sigma_{i,j} t$
    have no common nonzeros as well.
\end{itemize}
So, for some $i$ and $j$, 
the coloring $(g,\sigma_{i,j} t)$ has one less
essential arguments than
$(g,t)$.
If $(g,\sigma_{i,j} t)$ still has more than
$10$ essential arguments,
then (**) can be repeatedly applied to~$g$ and $\sigma_{i,j} t$,
and so on; 
finally,
we get two colorings~$g$, $\tilde t$
without common nonzeros
such that $(g,\tilde t)$ 
has exactly $10$ essential 
arguments and for some $i$, $j$
the coloring
$(g,\sigma_{i,j} \tilde t)$
has exactly~$11$ essential 
arguments, where $g$ and $\sigma_{i,j} \tilde t$
have no common nonzeros.
Denote by~$g_{10}$ and~$\tilde t_{10}$
the vertex colorings of $Q_{10}$
obtained from~$g$ 
and~$\tilde t$, respectively, by removing 
the nonessential arguments of $(g,\tilde t)$.
Since $g$ is a low-degree perfect coloring 
of $Q_n$, 
it is also true that $g_{10}$ 
is a low-degree perfect coloring of $Q_{10}$.
Similarly, $\tilde t_{10}$ is a fiber
of a low-degree perfect coloring of $Q_{10}$.
It follows that $(g_{10},\tilde t_{10})$,
up to equivalence, appeared
at step (ii)
of the algorithm when it ran
for $n=10$.
Moreover, from the consideration above
we know that $(g_{10},\tilde t_{10})$
is \emph{extendable} to the coloring 
$(g_{11},\tilde t_{11})$ of $Q_{11}$
by adding a nonessential argument so that
for some~$i$, the coloring $(g_{11},\sigma_{i,10}\tilde t_{11})$ has~$11$ essential arguments and the same number of colors
as $(g_{10},\tilde t_{10})$.
To confirm the claim of 
Proposition~\ref{p:n<11},
it is now sufficient
to check that all $k$-colorings
$(g,t)$ obtained
at step~(ii) of the algorithm
for $n=10$ are not extendable in the sense above.
This was checked 
during the run of the algorithm for $n=10$;
so, Proposition~\ref{p:n<11}
is confirmed.
\end{proof}

\begin{remark}
Alternatively, to confirm Proposition~\ref{p:n<11},
one can check that 
during the run of the algorithm for $n=11$,
every coloring $(g,t)$ that appears at step~(ii)
by combining a perfect coloring~$g$ and a fiber~$t$
has at least one nonessential argument. 
This needs more computing time,
but is also doable.
The most part of the proof above 
is still necessary to show that there is no
low-degree perfect coloring with $12$ or more essential arguments.
\end{remark}

\subsection
[Description of colorings into more than 2 colors]
{Description of colorings into more than $2$ colors}
\label{s:descr}

The following four propositions 
present constructions
of perfect colorings. 
The propositions
are straightforward; 
we skip their proofs.
The main result of this section,
Theorem~\ref{th:n-6+},
states that the class of low-degree 
perfect colorings is exhausted 
(up to equivalence and 
adding nonessential arguments)
by the colorings from Propositions~\ref{p:constr0}--\ref{p:constr3} and distance colorings of~$Q_2$
and~$Q_3$ with respect to one vertex.

\begin{proposition}\label{p:constr0}
Assume that we have a perfect
$k$-coloring~$f$ of $Q_{n-1}$ 
with quotient matrix~$S$.
Then the function~$g$,
$$
g(x_0,\ldots,x_{n-2},x_{n-1})
=
(x_{n-1},f(x_0,\ldots,x_{n-2}))
$$
is a perfect $2k$-coloring
of $Q_n$
with quotient matrix~$T$,
$
T = 
\begin{pmatrix}
 S  & \mathrm{Id} \\
 \mathrm{Id} & S
\end{pmatrix}
$.
Moreover, each $\lambda$-eigenvector
$\vc{v}$
of~$S$ corresponds to a 
$(\lambda+1)$-eigenvector 
$(\vc{v},\vc{v})$ and a 
$(\lambda-1)$-eigenvector 
$(\vc{v},-\vc{v})$ of~$T$.
In particular,
if the smallest eigenvalue
of~$S$ is
$\lambda_i(n-1,2)=(n-1)-2i$,
then 
the smallest eigenvalue
of~$T$ is $\lambda_{i+1}(n,2)=n-2(i+1)$.
\end{proposition}

\begin{proposition}\label{p:constr1}
Assume that each of $f$, $g$ is a
$\{0,1\}$-valued perfect
$2$- or $1$-colorings
of~$Q_{n-2}$
with quotient matrices,
respectively,  $\begin{pmatrix}
  n-2-b&b \\  b&n-2-b
\end{pmatrix}$ (or just~$(n-2-b)$, $b=0$ if $f$ is a $1$-coloring) and~$\begin{pmatrix}
     n-2-c&c   \\ c&n-2-c 
\end{pmatrix}$ (or~$(n-2-c)$, $c=0$, if $g$ is a $1$-coloring).
Then the function~$h$,
$$
h(x_0,\ldots,x_{n-3},x_{n-2},x_{n-1})
=
\begin{cases}
f(x_0,\ldots,x_{n-3}) &
\mbox{if $(x_{n-2},x_{n-1})=(0,0)$,}\\
1-f(x_0,\ldots,x_{n-3}) &
\mbox{if $(x_{n-2},x_{n-1})=(1,1)$,}\\
2+g(x_0,\ldots,x_{n-3}) &
\mbox{if $(x_{n-2},x_{n-1})=(0,1)$,}\\
3-g(x_0,\ldots,x_{n-3}) &
\mbox{if $(x_{n-2},x_{n-1})=(1,0)$,}
\end{cases}
$$
is a perfect $4$-coloring
of $Q_n$, 
with quotient matrix~$T$,
$$T = 
\left(\begin{array}{cccc}
 \cellcolor[gray]{0.95} n-2-b & \cellcolor[gray]{0.95} b & 1 & 1 \\
 \cellcolor[gray]{0.95} b & \cellcolor[gray]{0.95} n-2-b & 1 & 1 \\
 1 & 1 & \cellcolor[gray]{0.95} n-2-c & \cellcolor[gray]{0.95} c \\
 1 & 1 & \cellcolor[gray]{0.95} c & \cellcolor[gray]{0.95} n-2-c
\end{array}\right),
$$
which has eigenvalues $n$, $n-4$,
$n-2(b+1)$, $n-2(c+1)$
(with corresponding eigenvectors
$(1,1,1,1)$, $(1,1,-1,-1)$,
$(1,-1,0,0)$, $(0,0,1,-1)$).
\end{proposition}

To present the third construction, 
we define the following concept.
Two colors of 
a perfect $k$-coloring
are called \emph{twin}
if identifying them results in 
a perfect $(k-1)$-coloring
(of the same graph). For example,
in Proposition~\ref{p:constr1},
the first two colors of 
the coloring~$h$ are twin,
as well as the last two colors.
If we identify each of 
the two pairs of twin colors of~$h$,
then we get a perfect $2$-coloring
with quotient
matrix $\begin{pmatrix}
    n-2 & 2 \\ 2 & n-2
\end{pmatrix}$
and two essential arguments.
If $n\ge 4$, then there is another 
$2$-coloring with the same 
quotient matrix, but with 
four essential arguments. 
The third construction
is based on splitting 
each of the colors 
of that $2$-coloring into 
two twin colors.

We first define four perfect 
$3$-colorings $f_{a,b}$, $a,b\in\{0,1\}$, see Fig.~\ref{f:H52}.
\figureC
    Then, we define the following colorings:
    \begin{eqnarray}\label{eq:g4}
        g(x_0,\ldots,x_{n-1}) &=& f_{0,0}(x_0,\ldots,x_{3}),
        \quad n\ge 4,
    \\ \label{eq:gij}
     g_{i,j}(x_0,\ldots,x_{n-1}) &=& f_{x_i,x_j}(x_0,\ldots,x_{3}),
     \quad n\ge 5,\  i,j\in\{4,\ldots,n-1\}.
     \end{eqnarray}
\begin{proposition}\label{p:H42}
The colorings $g$ and $g_{i,j}$ (it is possible that $i=j$) of $Q_n$ are perfect with quotient matrices
$$
\left(\begin{array}{ccc}
    \cellcolor[gray]{0.95} n-3 & \cellcolor[gray]{0.95} 1 &  2 \\ 
    \cellcolor[gray]{0.95} 1 & \cellcolor[gray]{0.95} n-3 &  2 \\ 
     1 & 1 & \cellcolor[gray]{0.95} n-2
\end{array}\right),
\ n\ge 4, \qquad
\left(\begin{array}{ccc}
    \cellcolor[gray]{0.95} n-4 & \cellcolor[gray]{0.95} 2 & 2 \\ 
    \cellcolor[gray]{0.95} 2 & \cellcolor[gray]{0.95} n-4 & 2 \\ 
    1 & 1 & \cellcolor[gray]{0.95} n-2
\end{array}\right),
\ n\ge 5,
$$
with eigenvalues $n$, $n-4$, $n-4$ and $n$, $n-4$, $n-6$,
respectively.
The colorings $g$, $g_{i,i}$, and
$g_{i,j}$, $i\ne j$,
have $4$, $5$, and~$6$
essential arguments, respectively.
\end{proposition}
It is easy to see that the first two colors of~$g$ and~$g_{i,j}$ are twin.
Denote by~$f$ the perfect $2$-coloring obtained by identifying 
the these twin colors. It has $4$ essential arguments (see Fig.~\ref{f:F0}, the last diagram), and splitting the first color to produce~$g$ or~$g_{i,j}$ adds $0$, $1$, or $2$
essential arguments 
(respectively for $g$, $g_{i,i}$, and  $g_{i,i}$ with $i\ne j$).
Clearly, we can similarly split the last color of~$g$ or~$g_{i,j}$,
which can also add up to~$2$ essential arguments 
(say, $i'$ and~$j'$, which can be different from $i$, $j$).

\begin{proposition}\label{p:constr3}
By splitting each of the two colors of~$f$ in the way described above,
we get a perfect $4$-coloring 
of~$Q_n$ with 
the following quotient matrix:
$$
\left(\begin{array}{cccc}
    \cellcolor[gray]{0.95} n-b-2 & \cellcolor[gray]{0.95} b & 1 & 1 \\ 
    \cellcolor[gray]{0.95} b & \cellcolor[gray]{0.95} n-b-2 & 1 & 1 \\ 
    1 & 1  & \cellcolor[gray]{0.95} n-c-2 & \cellcolor[gray]{0.95} c \\
    1 & 1  & \cellcolor[gray]{0.95} c & \cellcolor[gray]{0.95} n-c-2 
\end{array}\right),
$$
where $b,c\in\{1,2\}$, $n\ge 3+\max(b,c)$, 
and the number of essential 
arguments is between 
$3+\max(b,c)$ and $2b+2c$.
The eigenvalues are 
$n$, $n-4$, $n-2b-2$,  $n-2c-2$.
\end{proposition}

\begin{theorem}[computational]\label{th:n-6+}
Every perfect coloring of $Q_n$ of degree at most~$3$
in more than~$2$ colors
is obtained by one 
of the constructions
in Propositions~\ref{p:constr0}--\ref{p:constr3}
or by adding inessential arguments
to one of the two perfect colorings
of $Q_2$ and $Q_3$, 
respectively,
with quotient matrices
\begin{equation}\label{eq:2mat}
\begin{pmatrix}
0 & 2 & 0 \\
1 & 0 & 1 \\
0 & 2 & 0
\end{pmatrix}
 \quad \mbox{and} \quad
\begin{pmatrix}
0 & 3 & 0 & 0 \\
1 & 0 & 2 & 0 \\
0 & 2 & 0 & 1 \\
0 & 0 & 3 & 0 
\end{pmatrix}
.
\end{equation}
\end{theorem}
Each coloring with quotient 
matrix from~\eqref{eq:2mat} 
is just the distance coloring
with respect to one vertex
(for example, the all-zero word).

\begin{theorem}[computational]\label{th:n-6}
Up to equivalence and 
adding inessential arguments,
there are $33$ perfect $2$-colorings
of hypercubes $Q_n$ with 
eigenvalue $n-6$; 
$20$ of them can be obtained by
identifying colors
in a perfect coloring with
more than~$2$ colors and 
minimum eigenvalue $n-6$
(described in Theorem~\ref{th:n-6+}).
The distribution by the number
of essential arguments is shown in Table~\ref{t:n-6}.
\end{theorem}
\begin{table}
$$
\begin{array}{r||c|c|c|c|c|c|c|c|c|c||c}
 n&1&2&3&4&5&6&7&8&9&10&\mbox{total} \\ \hline\hline
k=2 & 
0(1) & 0(1) & 1(1) & 0(1) & 1 &7 &12 & 9 & 2 & 1 & 33(4) \\ \hline
\scriptstyle 2' & \scriptstyle 0(1) & \scriptstyle 0(1) & \scriptstyle 1(1) & \scriptstyle 0(1) & \scriptstyle 1 &  \scriptstyle 6 & \scriptstyle 5 & \scriptstyle 5 & \scriptstyle 1 & \scriptstyle 1 & \scriptstyle 20(4) \\ \hline
k=3 &
0 & 0(1) & 0(1) & 1(1) & 1 & 2 & & & & & 4(3)
\\ \hline
k=4 &  
0 & 0(1) & 2(2) & 4(2) & 7 & 18 & 8 & 5 & 1 & 1 & 46(5)
\\ \hline
k=6 & 
0 & 0 & 1 & 1 & 1 & &&&&& 3
\\ \hline
k=8 & 
0 & 0 & 1 & 2 & 2 &&&&&& 5 \\ \hline\hline
k\ge 2 & 0(1) & 0(3) & 5(4) & 8(4) & 12 & 27 & 20 & 14 & 3 & 2 & 91(12)
\end{array}
$$
\caption{The number of perfect $k$-colorings in $Q_n$ with minimum eigenvalue $n-6$ (in the parenthesis, more than $n-6$) and without inessential arguments; the row ``$2'$'' indicates the numbers of $2$-colorings that can be subpartitioned into perfect colorings with more than two colors and eigenvalues not less than $n-6$.}
\label{t:n-6}
\end{table}

\section
[(n-4)-correlation-immune perfect colorings of Qn]
{$(n-4)$-correlation-immune perfect colorings of~$Q_n$}
\label{s:6-n}

By Lemma~\ref{l:deg-ci}(b),
$(n-4)$-resilient perfect $2$-colorings of $Q_n$ are in one-to-one correspondence 
with perfect $2$-colorings 
of degree at most~$3$, 
which are characterized 
in Section~\ref{s:k:n-6};
according to Table~\ref{t:n-6},
there are $1+3$, $1+4$, $2+4$, $9+4$, $21+4$, $30+4$, $32+4$, and $33+4$ equivalence classes of such $2$-colorings for $n=3$, $4$, \ldots, $9$, and $n\ge10$, respectively, where $a+b$ means that $a$ colorings are $(n-4)$-resilient but not 
 $(n-3)$-resilient (and hence have quotient matrix $\quot{3&n-3\\n-3&3}$) while $b$ colorings are at least $(n-3)$-resilient.

It remains to consider the case of $2$-colorings with non-symmetric quotient matrix 
and the case of colorings with more than~$3$ colors. By Corollaries~\ref{c:FDF} and~\ref{c:FDF2}, in both cases we have $n\le 9$.
Each such a coloring has 
an $(n-4)$-correlation-immune fiber
with density of ones strictly 
between~$0$ and~$1/2$.
Such Boolean functions are not resilient 
and they are not included in the classification~\cite{Rasoolzadeh2023}.
However, $t$-correlation-immune 
Boolean functions are well known to be equivalent to simple binary
orthogonal arrays of strength~$t$
(essentially, such an array 
can be defined as the set of $1$s 
of a $t$-correlation-immune Boolean function),
and the non-resilient $(n-4)$-correlation-immune Boolean functions can be extracted 
(see Table~\ref{t:ci}) 
\begin{table}[t]
\mbox{}\hfill
$
\begin{array}{r||cc|cc|cc|cc|cc|cc|cc}
\mbox{density:} 
&  \makebox[1mm][l]{1/16} 
&& \makebox[1mm][l]{1/8 } 
&& \makebox[1mm][l]{3/16} 
&& \makebox[1mm][l]{1/4 } 
&& \makebox[1mm][l]{5/16} 
&& \makebox[1mm][l]{3/8 } 
&& \makebox[1mm][l]{7/16} 
\\ \hline
n = 9: & - && - && - && 2 &(2)& - && - && - & \\
n = 8: & - && - && - && 3 &(3)& - && - && - & \\
n = 7: & - && 1&(1) & 1&(0) & 13&(7) & 7&(0) & 82&(10) & 71 &(0)  \\
n = 6: & - && 1 &(1) & 2 &(0) & 21 &(9) & 22 &(0) & 178 &(6) & 261 &(0) \\
n = 5: & 1 & (1) & 4 & (4) & 8 &(1) & 36 & (9) & 73 &(3) & 190 &(3) & 266 &(0) \\\hline
\end{array}
$
\hfill\mbox{}
\caption{The number of equivalence classes 
of $(n-4)$-correlation-immune Boolean function
with density of $1$s less than $1/2$
(in parenthesis, the number of functions that are fibers of $(n-4)$-correlation-immune perfect colorings).
}
\label{t:ci}
\end{table}
from 
the classification 
of orthogonal arrays of small parameters by Schoen, Eendebak, and Nguyen~\cite{SEN:2010:OA}.
With that classification, for each~$n$, 
$n\le 9$, we can apply the algorithm similar 
to Section~\ref{s:class}
to find all equitable partitions
with correlation immunity~$n-4$.

In Sections~\ref{s:H9}--\ref{s:H5},
we describe the obtained classification 
of $(n-4)$-correlation-immune perfect colorings of~$Q_n$ for $n=5$, $n=6$, $n=7$, $n=8$, and $n=9$.
The colorings with higher~$n$ have also higher correlation immunity order~$n-4$, and we can say more about their structure. 
We pay a special attention to perfect 
colorings with only two eigenvalues~$n$ and~$n-6$, which exist only for $n\in\{5,7,9\}$,
apart of the resilient $2$-colorings
(it is not difficult to conclude from \eqref{eq:rhoij} that non-resilient $2$-colorings with eigenvalue~$n-6$ do not exist
for $n=6$ and $n=8$; 
by Corollary~\ref{l:union}, this also means, for $n\in \{6,8\}$,
the nonexistence of perfect colorings with 
more than two colors and only two eigenvalues~$n$ and~$n-6$).

\subsection
[Perfect colorings of Q5 with correlation immunity 1]
{Perfect colorings of $Q_5$ with correlation immunity $1$}\label{s:H5}

There are many $1$-correlation-immune 
perfect colorings of~$Q_5$, and in the following theorem,
we reflect only the number of such $k$-colorings for each~$k$. 
After that, we will focus our attention
on the special case of only 
two eigenvalues~$5$ and~$1$.

\begin{theorem}[computational]
    The number of equivalence classes of $1$-cor\-re\-la\-tion-immune perfect $k$-colorings
    of $Q_5$ for $k=2,3,\ldots,16$
    is, respectively,
    $7(6)$, $9$, $31(22)$, $14$, $24$, $14$, $17(16)$, $1$, $2$, $0$, $1$, $0$, $0$, $0$, $1(1)$ (in the parenthesis, the number of
    only resilient colorings is given);
    the number of different quotient matrices
    is, respectively, 
    $5(4)$, $7$, $19(10)$, $10$, $16$, $10$, $12(11)$, $1$, $2$, $0$, $1$, $0$, $0$, $0$, $1(1)$.
\end{theorem}

The only (up to reordering of colors) non-symmetric quotient matrix 
of a perfect $2$-coloring 
of $Q_5$ with 
eigenvalue~$1$ is
$\begin{pmatrix}
    2&3\\1&4
\end{pmatrix}$. 
We see that there are $8$ vertices
of first color, and they induce a cycle 
or the union of two cycles 
in~$Q_5$.
By Lemma~\ref{l:antipod},
it can be 
only the union of two mutually antipodal
$4$-cycles. 
Hence, a coloring is unique up to equivalence. 
By Corollary~\ref{l:union},
any perfect $3$- or $4$-coloring of~$Q_5$
with two eigenvalues~$5$ and~$1$ have,
respectively, $2$ or~$4$ fibers 
that are perfect $2$-colorings considered above. 
Such perfect $3$- or $4$-colorings can be easily
characterized; we skip the proof, noting 
also that they have degree~$2$ and hence are
included in the characterization in Section~\ref{s:k:n-6}.

\figureB
 \begin{theorem}
 A perfect coloring of $Q_5$ with eigenvalues
 $5$ and~$1$ has at least one nonessential argument
 and one of the following quotient matrices, up to 
 ordering of colors:
 $\begin{pmatrix}2&3\\1&4\end{pmatrix}
 $,
 \quad
 $
 \begin{pmatrix}3&2\\2&3\end{pmatrix},
 \quad
 \quot{2&1&2\\1&2&2\\1&1&3}
 $,
 \quad
 $
 \quot{2&1&1&1\\1&2&1&1\\
                1&1&2&1\\1&1&1&2}
 $.
 The number of equivalence classes of such colorings is, respectively, $1$, $2$, 
 $2$,~$3$.
 Representatives of the $3$ equivalence classes of $4$-colorings are obtained from the colorings
 of $Q_4$ at Fig.~\ref{f:H42} by adding 
 a nonessential argument; representatives
 of $3$- and $2$-colorings are obtained 
 from those three $4$-colorings by merging colors.
 \end{theorem}

\subsection
[Perfect colorings of Q6 with correlation immunity 2]
{Perfect colorings of $Q_6$ with correlation immunity $2$}
\label{s:H6}

The results 
of the characterization $2$-correlation-immune perfect colorings of~$Q_6$
are indicated in Table~\ref{tab:Q6}.
\begin{table}[t]
    \centering
\newcommand\qwe{\!\!\!}
\newcommand\qw{\!\!}
$\quot{
0&0&1&1&1&1&1&1\\
0&0&1&1&1&1&1&1\\
1&1&0&0&1&1&1&1\\
1&1&0&0&1&1&1&1\\
1&1&1&1&0&0&1&1\\
1&1&1&1&0&0&1&1\\
1&1&1&1&1&1&0&0\\
1&1&1&1&1&1&0&0
}^{\qwe7}$
$\quot{
0&0&1&1&1&1&2\\
0&0&1&1&1&1&2\\
1&1&0&0&1&1&2\\
1&1&0&0&1&1&2\\
1&1&1&1&0&0&2\\
1&1&1&1&0&0&2\\
1&1&1&1&1&1&0
}^{\qwe5}$
$\quot{
0&0&1&1&2&2\\
0&0&1&1&2&2\\
1&1&0&0&2&2\\
1&1&0&0&2&2\\
1&1&1&1&0&2\\
1&1&1&1&2&0
}^{\qwe6}$
$\quot{
0&0&1&1&2&2\\
0&0&1&1&2&2\\
1&1&0&0&2&2\\
1&1&0&0&2&2\\
1&1&1&1&1&1\\
1&1&1&1&1&1
}^{\qwe15}$
$\quot{
0&0&1&1&4\\
0&0&1&1&4\\
1&1&0&0&4\\
1&1&0&0&4\\
1&1&1&1&2
}^{\qwe4}$
$\quot{
0&0&2&2&2\\
0&0&2&2&2\\
1&1&0&2&2\\
1&1&2&0&2\\
1&1&2&2&0
}^{\qwe2}$
$\quot{
0&0&2&2&2\\
0&0&2&2&2\\
1&1&0&2&2\\
1&1&2&1&1\\
1&1&2&1&1
}^{\qwe9}$
$\quot{
0&0&2&4\\
0&0&2&4\\
1&1&0&4\\
1&1&2&2
}^{\qwe2}$
$\quot{
0&0&3&3\\
0&0&3&3\\
1&1&2&2\\
1&1&2&2
}^{\qwe9}$
${\boldmath\quot{
0&2&2&2\\
2&0&2&2\\
2&2&0&2\\
2&2&2&0
}}^{\qwe2}$
$\quot{
0&2&2&2\\
2&0&2&2\\
2&2&1&1\\
2&2&1&1
}^{\qwe6}$
$\quot{
1&1&2&2\\
1&1&2&2\\
2&2&1&1\\
2&2&1&1
}^{\qwe19}$
$\quot{
0&0&6\\
0&0&6\\
1&1&4
}^{\qw1}$
${\boldmath\quot{
0&2&4\\
2&0&4\\
2&2&2
}}^{\qw2}$
$\quot{
0&3&3\\
2&2&2\\
2&2&2
}^{\qw5}$
$\quot{
1&1&4\\
1&1&4\\
2&2&2
}^{\qw4}$
${\boldmath\QQuot{
0&\!\!\!6\\[-0.6ex]
2&\!\!\!4
}}^{1}$
${\boldmath\QQuot{
2&\!\!\!4\\[-0.6ex]
4&\!\!\!2
}}^{2}$
$\QQuot{
3&\!\!\!3\\[-0.6ex]
3&\!\!\!3
}^{9}$
\\
${\boldmath\QQuot{
0&\!\!\!6\\[-0.6ex]
6&\!\!\!0
}}^{1}$
${\boldmath\QQuot{
1&\!\!\!5\\[-0.6ex]
5&\!\!\!1
}}^{1}$
${\boldmath\QQuot{
1&\!\!\!5\\[-0.6ex]
3&\!\!\!3
}}^{1}$
${\boldmath\quot{
0&3&3\\
2&1&3\\
2&3&1
}}^{\qw1}$
$\quot{
0&0&6\\
0&0&6\\
3&3&0
}^{\qw2}$
$\quot{
0&0&3&3\\
0&0&3&3\\
1&1&1&3\\
1&1&3&1
}^{\qw1}$
$\quot{
0&0&3&3\\
0&0&3&3\\
3&3&0&0\\
3&3&0&0
}^{\qw9}$
$\quot{
0&1&2&3\\
1&0&3&2\\
2&3&0&1\\
3&2&1&0
}^{\qw2}$
    \caption{Quotient matrices of $2$-correlation-immune perfect colorings of $Q_6$ with the number on nonequivalent colorings (indicated as the degree); the matrices marked by bold have better correlation immunity.}
    \label{tab:Q6}
\end{table}
The remarkable group of matrices 
in the first three rows of Table~\ref{tab:Q6}
is characterized by the following property:
each of these quotient matrices can be obtained by merging some groups of colors of a perfect $8$-coloring with the first quotient matrix
(an example of such a $8$-coloring
is $f(x_1,\ldots,x_8) = (x_1+x_2+x_4+x_5,x_2+x_3+x_5+x_6,x_1+x_2+x_3)$).
However, 
similar to Section~\ref{s:H7b}, 
not all colorings with these matrices 
have a perfect $8$-coloring refinement.

\subsection
[Perfect colorings of Q7 with correlation immunity 3]
{Perfect colorings of $Q_7$ with correlation immunity~$3$}\label{s:H7}

The $3$-correlation-immune perfect colorings 
of $Q_7$ are divided into three groups:
\begin{itemize}
    \item perfect colorings of correlation immunity at least~$4$ with quotient matrices
    $\Quot{0&7\\7&0}$ ($1$ equivalence class),
    $\Quot{1&6\\6&1}$ ($1$ equivalence class),
    $\Quot{2&5\\5&2}$ ($2$ equivalence classes);
    \item perfect colorings with exactly two eigenvalues~$7$ and~$-1$; we consider this class of colorings in separate Subsection~\ref{s:H7b} below;
    \item perfect colorings with more than two eigenvalues, described in the following theorem.
\end{itemize}

\begin{theorem}
Up to equivalence, the perfect colorings
of $Q_7$ with correlation immunity~$3$
and more than~$2$ eigenvalues are exhausted by the following list:
\begin{itemize}
    \item 
    Colorings with the following quotient matrices
(respectively, $2$, $9$, $1$, $1$ colorings):
$\quot{0&1&6\\1&0&6\\3&3&1}
,\quad
\quot{0&1&3&3\\1&0&3&3\\3&3&0&1\\3&3&1&0}
,\quad
\quot{1&3&3\\2&1&4\\2&4&1}
,\quad
\quot{0&1&3&3\\1&0&3&3\\1&1&1&4\\1&1&4&1}.$ 
\item
Two colorings with quotient matrix
$\quot{0&3&4\\3&0&4\\2&2&3}$ and their refinements with the following matrices
(respectively, $3$, $7$, $2$, $7$, and $5$ colorings):
$$\quot{0&3&2&2\\3&0&2&2\\2&2&0&3\\2&2&3&0}
,\quad
\quot{0&3&2&2\\3&0&2&2\\2&2&1&2\\2&2&2&1}
,\quad
\quot{0&3&1&3\\3&0&1&3\\2&2&0&3\\2&2&1&2}
,\quad
\quot{0&3&1&1&2\\3&0&1&1&2\\2&2&0&1&2\\2&2&1&0&2\\2&2&1&1&1}
,\quad
\quot{0&3&1&1&1&1\\3&0&1&1&1&1\\2&2&0&1&1&1\\2&2&1&0&1&1\\2&2&1&1&0&1\\2&2&1&1&1&0}.$$
\end{itemize}
\end{theorem}
It is notable that for a perfect coloring with 
one of the last four quotient matrices, merging the first two colors results in 
a perfect coloring from the following subsection.
\subsubsection
[Perfect colorings of Q7 with two eigenvalues 7 and -1]
{Perfect colorings of $Q_7$ with two eigenvalues $7$ and $-1$}\label{s:H7b}
Perfect $2$-colorings of $Q_7$ with 
eigenvalue~$-1$, i.e., with quotient matrix
$\begin{pmatrix}
    \mu-1 &  7-\mu+1\\  \mu& 7-\mu
\end{pmatrix}$
are essentially special cases 
of $\mu$-fold $1$-perfect codes.
A~set~$C$ of vertices of an $r$-regular graph
is called a \emph{$\mu$-fold $1$-perfect code}
(\emph{multifold $1$-perfect code};
in the case $\mu=1$, just
\emph{$1$-perfect code})
if every radius-$1$ ball contains exactly
$\mu$ elements from~$C$.
Equivalently, 
the characteristic function of~$C$ is
a perfect coloring with quotient matrix
$\begin{pmatrix}
     \mu-1 &  r-\mu+1\\  \mu& r-\mu
\end{pmatrix}$.
A $1$-perfect codes in $Q_7$,
the Hamming code, is unique up to equivalence
and there are exactly $240$ such 
equivalent codes.
Two-fold $1$-perfect codes 
in $Q_7$ 
were essentially characterized 
in~\cite[Section~4]{Kro:small:2017}
(the case when the code can be split into two 
$1$-perfect codes)
and~\cite[Section~VI.D]{KroPot:multifold}
(the unsplittable case).
From those results and 
our independent calculation,
we know that there are exactly three 
two-fold $1$-perfect codes 
in $Q_7$ (below, each code is represented by the list of values of its characteristic Boolean function, packed in the hexadecimal form, 
$\mathtt{0}=0000$, $\mathtt{1}=0001$, \ldots,
$\mathtt{f}=1111$, and followed by the automorphism group order):

$\mathtt{c30000c3003c3c00003c3c00c30000c3}$, $|\mathrm{Aut}|=1536$;

$\mathtt{c30000a5005a3c00003c5a00a50000c3}$, $|\mathrm{Aut}|=128$;

$\mathtt{c2100426021c91800189384064200843}$, $|\mathrm{Aut}|=96$.

A $3$-fold $1$-perfect code,
as we see from the quotient matrix,
induces a $2$-factor
(the union of disjoint cycles)
in the ambient graph.
We found in $Q_7$ that the number 
of equivalence classes 
of such codes is~$9$. 
Below, we list representatives
together with the orders 
of automorphism groups 
and cycle lengths 
(e.g., $4^2 10^4$ denotes $2$ cycles of 
length~$4$ and $4$ cycles of length~$10$;
the self-complementary,
in view of Lemma~\ref{l:antipod},
cycles are indicated by $\overline{14}$):

$\mathtt{e70000db00bd7e00007ebd00db0000e7}$, $|\mathrm{Aut}|=384$, cycles: $6^8$;

$\mathtt{e610049b029d7680016eb940d9200867}$, $|\mathrm{Aut}|=32$, cycles: $12^4$;

$\mathtt{e610049b025eb58001ad7a40d9200867}$, $|\mathrm{Aut}|=8$, cycles: $24^2$;

$\mathtt{e424128709395ac0035a9c90e1482427}$, $|\mathrm{Aut}|=12$, cycles: $24^2$;

$\mathtt{f00c0c3303c33cc0033cc3c0cc30300f}$, $|\mathrm{Aut}|=480$, cycles: $4^{12}$;

$\mathtt{f009096906996660066699609690900f}$, $|\mathrm{Aut}|=160$, cycles: $4^2 10^4$;

$\mathtt{f00c0c3303a55ac0035aa5c0cc30300f}$, $|\mathrm{Aut}|=16$, cycles: $4^4 16^2$;

$\mathtt{f00c0963063699c003996c60c690300f}$, $|\mathrm{Aut}|=16$, cycles: $4^2 10^4$;

$\mathtt{e428141b03659ac00359a6c0d8281427}$, $|\mathrm{Aut}|=8$, cycles: $\overline{14}\vphantom{0}^2 10^2$.

The first $4$ of the codes above can be split
into three $1$-perfect codes,
while the last~$5$ cannot, 
as can be seen from the cycle lengths
(an example of such unsplittable $3$-fold
$1$-perfect code in~$Q_7$ was found 
in~\cite{Westerback:2007}; in~$Q_{2^m-1}$, 
$m\ge 4$, unsplittable multifold
$1$-perfect codes were constructed 
in~\cite{KroPot:nonsplittable}).

Since the complement of 
a $\mu$-fold $1$-perfect code in~$Q_7$ 
is a $(8-\mu)$-fold $1$-perfect code
and the $4$-fold case 
is characterised above 
(by Lemma~\ref{l:deg-ci}(b), $4$-fold $1$-perfect codes in~$Q_7$ are equivalent to perfect colorings with quotient matrix $\begin{pmatrix}
    4&3\\3&4
\end{pmatrix}$; from Table~\ref{t:n-6} we see that there are $21=1+1+7+12$ equivalence classes of such colorings),
the classification of multifold
$1$-perfect codes, as well as perfect 
$2$-colorings with eigenvalue~$-1$,
in~$Q_7$ is completed.
Based on the classification, 
we can state the following.
\begin{proposition}[computational]\label{p:split}
In $Q_7$, there are 
$3$-fold, $4$-fold, and $5$-fold
$1$-perfect codes that cannot be split 
into $1$-perfect codes;
all $2$-fold and $6$-fold 
$1$-perfect codes are splittable 
in this sense.
Each $\mu$-fold $1$-perfect code 
in $Q_7$, $\mu\ne 3$,
includes as a subset at least one 
$1$-perfect code
(i.e., can be split into $1$-perfect code
and $(\mu-1)$-fold $1$-perfect code).
\end{proposition}

Based on Proposition~\ref{p:split},
we can now classify all partitions
of the vertex set 
of~$Q_7$ into
multifold $1$-perfect codes,
or, equivalently,
perfect coloring of~$Q_7$
with only two eigenvalues~$-1$ and~$7$.
We will call such partition
into $\mu_0$-fold, $\mu_1$-fold, 
\ldots, $\mu_{k-1}$-fold $1$-perfect codes
a \emph{$(\mu_0,\ldots,\mu_{k-1})$-partition}.
For example, we know that every 
$(3,3,1,1)$-partition can be obtained
from some $(4,4)$-partition $(C_0,C_1)$
by splitting each of $C_0$, $C_1$ into
a $1$-perfect code and a $3$-fold $1$-perfect
code 
(which is computationally easy 
taking into account that there are 
only~$240$ $1$-perfect codes).
On the other hand, every 
$(3,3,2)$-partition can be obtained
from some $(3,3,1,1)$-partition 
$(C_0,C_1,C_2,C_3)$
by unifying~$C_2$ and~$C_3$.

\begin{theorem}
[computational]
\label{p:-1}
The total numbers 
(up to reordering of codes /
renaming the colors)
of the partitions of the vertex set 
of~$Q_7$ into
multifold $1$-perfect codes,
corresponding to perfect colorings
with eigenvalues~$7$ and~$-1$
are the following 
(scaled cardinality spectrum : 
number of equivalence classes : 
total number):\\
$$
\begin{array}{lll}
(7,1):1:240,         &   (4,2,2):40:2029230,      &   (4,1,1,1,1):26:771120,      \\
(6,2):3:12180,       &   (3,3,2):69:5444880,      &   (3,2,1,1,1):108:9770880,    \\
(5,3):9:322896,      &   (5,1,1,1):6:181440,      &   (2,2,2,1,1):129:8342040,    \\
(4,4):21:505715,     &   (4,2,1,1):58:4331880,    &   (3,1,1,1,1,1):27:1069824,   \\
(6,1,1):3:12600,     &   (3,3,1,1):81:5866560,    &   (2,2,1,1,1,1):99:4621680,   \\
(5,2,1):9:519120,    &   (3,2,2,1):121:13401360,  &   (2,1,1,1,1,1,1):29:685440,  \\
(4,3,1):38:4537680,  &   (2,2,2,2):66:1890105,    &   (1,1,1,1,1,1,1,1):11:27360. 
\end{array}
$$
\end{theorem}
Representatives of equivalence classes
can be found in~\cite{Perfect-related}.

\subsection
[Perfect colorings of Q8 with correlation immunity 4]
{Perfect colorings of $Q_8$ with correlation immunity~$4$}\label{s:H8}

The case $n=8$ is special.
According to Table~\ref{t:ci},
any $4$-correlation-immune coloring
of~$Q_8$ that is not a resilient $1$- 
or $2$-coloring
has a fiber of density~$1/4$.
A $4$-correlation-immune Boolean function 
in $n$ arguments with density~$1/4$ belongs
to a special class, for which it was proved
in \cite[Theorem~1]{Kro:OA13}
that such a function is necessarily
a fiber of a perfect $3$-coloring.
This explains why the variety 
of perfect colorings 
described in the following classification theorem 
is so small.

\begin{theorem}[computational]
Up to equivalence, the perfect colorings
of $Q_8$ with correlation immunity~$4$,
apart of the $4$-resilient $2$-colorings,
are exhausted by three colorings with quotient matrix
$\quot{0&2&6\\2&0&6\\3&3&2}$ and $11$ 
their refinements with quotient matrix 
$\quot{0&2&3&3\\2&0&3&3\\3&3&0&2\\3&3&2&0}$.
\end{theorem}

\subsection
[Perfect colorings of Q9 with correlation immunity 5]
{Perfect colorings of $Q_9$ with correlation immunity~$5$}\label{s:H9}
According to \cite{FDF:CorrImmBound},
every unbalanced Boolean function in $n$ arguments
of correlation immunity $2n/3-1$
is a perfect $2$-coloring.
Such functions with $n=9$ were characterized 
by Kirienko~\cite{Kirienko2002}; it was found that
there are exactly~$2$ equivalence classes of them.
In particular, there are exactly~$2$ equivalence classes of equitable $2$-partitions with quotient matrix $\begin{pmatrix}
    0&9\\3&6
\end{pmatrix}$.

There are representatives of those 
equivalence classes
with an easy group description:
the vertices $(x_0,\ldots,x_8)$
of the first color are the solutions
of the equation
$$ (x_0+x_1, x_0+x_2)\star
   (x_3+x_4, x_3+x_5) =
   (x_6+x_7, x_6+x_8),  $$
where $+$ is the addition over $\mathbb Z_2$ 
and $( \{(0,0),(0,1),(1,0),(1,1)\},\star )$
is a group isomorphic 
to~$\mathbb Z_2^2$ 
or~$\mathbb Z_4$.

By Corollaries~\ref{c:FDF} and~\ref{l:union},
the only admissible quotient matrices for
$5$-correlation-immune perfect
colorings of $Q_9$ with 
more than two colors are
\begin{equation}\label{eq:34}
  \begin{pmatrix}
   0&3&6\\3&0&6\\3&3&3 
\end{pmatrix}
\quad \mbox{and} \quad
\begin{pmatrix}
   0&3&3&3\\3&0&3&3\\3&3&0&3\\3&3&3&0 
\end{pmatrix},  
\end{equation}
up to reordering of colors.
Apart of the general approach mentioned in the beginning of Section~\ref{s:6-n},
there is a rather simple
algorithm of obtaining 
such $3$-colorings from $32$ inequivalent 
perfect $2$-colorings with quotient matrix
$
\begin{pmatrix}
   3&6\\6&3 
\end{pmatrix}
$, 
because we need 
to subpartition
one of the colors into two independent sets.
To do this,
we consider the subgraph of $Q_9$
induced by the vertices of one of the two colors
and the partitions 
of those vertices into bipartite
sets of that subgraph 
(which can be disconnected,
so such a partition is not unique in general).
Not all such partitions result 
in a perfect $3$-coloring; 
in particular, we are interested only in
self-complimentary partitions 
(see Lemma~\ref{l:antipod}),
which simplifies computations.
Similarly, from perfect $3$-colorings, we can obtain all nonequivalent perfect $4$-colorings
with required quotient matrix.
\begin{theorem}[computational]\label{p:934}
Up to equivalence, there are $11$
perfect $3$-colorings
and $10$ 
perfect $4$-colorings of $Q_9$
with quotient matrices~\eqref{eq:34}.
Each of the $3$-colorings can be obtained from 
a perfect $4$-coloring by unifying two colors.

  Nine of the $10$ perfect $4$-colorings can 
be obtained from the $9$ inequivalent perfect $2$-colorings~$g$ of $Q_6$
with quotient matrix
$\big(\begin{smallmatrix}
    3&3 \\ 3&3
\end{smallmatrix}\big)$
as
$$
f(x_0,x_1,x_2,y_0,\ldots,y_5) =
(g(y_0,\ldots,y_5)+x_0+x_1+x_2,y_0+\ldots+y_5),
$$
where $+$ is the modulo-$2$ addition.
The remaining perfect $4$-coloring has the form 
$$
f(x_0,\ldots,x_9) = 
   \star((x_0+x_1, x_0+x_2),
   (x_3+x_4, x_3+x_5),
   (x_6+x_7, x_6+x_8)),
   $$
where $+$ is the addition over $\mathbb Z_2$ 
and $\star(\cdot,\cdot,\cdot)$
is a ternary quasigroup operation 
on $\{(0,0),(0,1),(1,0),(1,1)\}$,
the composition 
$\star(a,b,c) = (a * b)\circ c $
of two different groups isomorphic 
to~$\mathbb Z_4$.
\end{theorem}
Representatives of perfect colorings 
described in Theorem~\ref{p:934}
can be found in~\cite{Perfect-related}. 
More information about perfect colorings
of $Q_{3m}$ with eigen\-val\-ue~$-m$ 
and corresponding resilient 
Boolean functions can be found 
in~\cite{Krotov:ISIT2019:Resilient}.

\section{Conclusion}
In this paper, we classified two families of perfect colorings of the hypercube~$Q_n$, perfect colorings of degree at most~$3$ and perfect colorings of correlation immunity (at least) $n-4$.
Perfect $2$-colorings of degree at most~$3$ are in one-to-one correspondence (Lemma~\ref{l:deg-ci}(b)) 
with $(n-4)$-resilient perfect $2$-colorings. While the class
of such functions is characterized computationally,
we could not find a theoretical construction for all of them. However, even without such construction, 
we can conclude that $37$ of $961$ Boolean
functions of degree at least~$3$ in $10$ 
or more arguments 
have the very regular structure 
of a perfect $2$-coloring 
(additionally, $1$ function is a perfect $1$-coloring),
including the unique such function with $10$
essential arguments.
The same can be said about $38$ of $961$
$(n-4)$-resilient Boolean functions.

The situation with perfect colorings of~$Q_n$ into more than two colors is different. There is no direct connection between such colorings
of degree at most~$3$ and of correlation immunity at least $n-4$, and we comment them separately. The class of such colorings
of degree at most~$3$ is infinite, because
adding a nonessential argument does not
change the degree or the perfectness of the coloring. However, the number of essential
arguments is limited by~$10$, and hence
the number of perfect colorings with all
essential arguments is finite; this number
is $66$ (for all $n$ and any number of colors larger than~$2$). We found theoretical constructions that describe all 
these colorings. 
It should be noted that non-perfect degree-$3$ colorings
with more than $2$ colors (where each color is represented by a Boolean function of degree at most~$3$) are not described 
yet, and we do not even know if the number
of essential arguments of such a coloring
is bounded by~$10$.

In contrast to the degree-$3$ 
perfect colorings, 
$(n-4)$-correlation-immune perfect colorings
with more than~$2$ colors cannot be constructed by adding nonessential arguments,
and we classified them separately for
each $n$ from~$5$ to~$9$, 
paying more attention for higher~$n$
because the strength $n-4$ is also higher 
and the colorings are more structured
(in particular, any non-resilient $(n-4)$-correlation-immune Boolean function in $n=8$ or $n=9$ arguments is necessarily a fiber of a perfect coloring).
It should be noted that we can always
treat colors as binary tuples (for example,
for $4$-colorings, the colors can be $00$, $01$, $10$, and $11$), and colorings
into more than~$2$ colors can be treated as 
vectorial Boolean functions, in the spirit of~\cite{Carlet:vect}.
In that sense, perfect colorings 
with symmetric quotient matrices
can be treated as resilient vectorial Boolean functions
(as follows from~\eqref{eq:rhoij}, all densities are equal if the matrix is symmetric). In Section~\ref{s:6-n}, 
we see several classes of such 
$(n-4)$-resilient $4$- and $8$-colorings with different quotient matrices.

\subsection*{Acknowledgments}
The authors thank Pieter Eendebak for sharing the database with orthogonal arrays.
The study was funded by the Russian Science Foundation, grant 22-11-00266, \url{https://rscf.ru/en/project/22-11-00266/}.

\subsection*{Declaration of competing interest}
The authors declare that they have no known competing financial interests or personal
relationships that could have appeared to influence the work reported in this paper.

\subsection*{Data availability}
In our computations, we
used the list of the $761$ nonequivalent
$6$-resilient Boolean functions 
in $10$ arguments,
generated by Rasoolzadeh \cite{Rasoolzadeh2023}
and available at
\url{https://gitlab.science.ru.nl/shahramr/ResilientFunctions.git},
file \texttt{F10T6.txt}. 
The dataset containing the results 
of the classifications described in 
Sections~\ref{s:k:n-6} and~\ref{s:6-n} 
is available in the IEEE DataPort
repository~\cite{Perfect-related}.


\providecommand\href[2]{#2} \providecommand\url[1]{\href{#1}{#1}}
  \def\DOI#1{{\href{https://doi.org/#1}{https://doi.org/#1}}}\def\DOIURL#1#2{{\href{https://doi.org/#2}{https://doi.org/#1}}}

\end{document}